\documentclass[10pt]{article}
\usepackage{amsfonts, latexsym, amscd, graphics, amssymb}

\makeatletter
\@addtoreset{equation}{section}
\makeatother

\def\text#1{\mbox{\rm #1}}
\newcommand{\be}{\begin{equation}}
\newcommand{\ee}{\end{equation}}
\newcommand{\ba}{\begin{eqnarray}}
\newcommand{\ea}{\end{eqnarray}}

\newcommand{\pa}{\partial}
\newcommand{\f}{\frac}

\begin{document}

\title{ Birkhoffian formulation of the dynamics of LC circuits}
\date{}
\date{}
\author{\large{Delia Ionescu\footnote{{\it Permanent address}: Institute of Mathematics,
Romanian Academy of Sciences,
 P.O. Box 1-764, RO-014700, Bucharest, Romania, Delia.Ionescu@imar.ro} , J\"{u}rgen Scheurle}\footnote{, * Research supported by
the EU through the Research Training Network \textit{Mechanics and
Symmetry in Europe}(MASIE)}
\\
{\small {\it Zentrum Mathematik      der Technische Universit\"{a}t M\"{u}nchen,}}\\
                   {\small {\it D-85747 Garching bei M\"{u}nchen,
                   Germany}}}
\maketitle
\maketitle

\textit{\textbf {Abstract.}} We present a formulation of general
nonlinear LC circuits within the framework of Birkhoffian
dynamical systems on manifolds. We develop a systematic procedure
which allows, under rather mild non-degeneracy conditions,  to
write the governing equations for the mathematical description of
the dynamics of an LC circuit as a Birkhoffian differential
system. In order to illustrate the advantages of this approach
compared to known Lagrangian or Hamiltonian approaches we discuss
a number of specific examples. In particular, the Birkhoffian
approach includes networks which contain closed loops formed by
capacitors, as well as inductor cutsets. We also extend our
approach to the case of networks which contain independent voltage
sources as well as independent current sources. Also, we derive a
general balance law for an associated "energy function".


\section{Introduction}

In this paper we give a formulation of the dynamics of LC circuits
within the framework of  Birkhoffian systems \cite{bir}. Based on
the constitutive relations of the involved inductors and
capacitors and on  Kirchhoff's laws, we  define  a configuration
space and a corresponding Birkhoffian that describes the
``elementary work'' done by a set of  ``generalized forces''. As a
matter of fact, in
order to cover circuits for which the topological assumptions usually imposed in the literature, are not satisfied,  we are forced 
to describe a single circuit by a whole family of Birkhoffian systems parameterized by a finite number of real parameters. Relevant 
values of these parameters correspond to initial values for the
time evolution of certain state variables of the circuit. The
dimension of each configuration space is given by the cardinality
of a selection of loops that cover the whole circuit.

In order to study  the dynamics of LC circuits, various Lagrangian and Hamiltonian formulations have been considered in the 
literature (see for example \cite{bern}, \cite{bloch},  \cite{moser}, \cite{chua}, \cite{maschke}, \cite{maschke2}).

\noindent In the Lagrangian approach, a central issue is the selection of suitable coordinates and corresponding velocities in terms of which 
the Lagrangian function is expressed.
A specific technique for the sometimes difficult task of  choosing the proper Lagrangian variables is presented in paper 
\cite{chua}.

\noindent More often Hamiltonian formulations have been used to describe circuit equations. In \cite{bern} it is shown how to construct, based 
on the circuit topology,  canonical variables  and a  Hamiltonian, so that the circuit equations attain canonical form.

\noindent For a more 
general approach including also resistors, the RLC circuits, see Brayton-Moser's
approach \cite{moser}. In \cite{moser}, under the hypothesis that
the currents through the inductors and the voltages across the
capacitors determine all currents and voltages in the circuit via
Kirchhoff's law, is proved the existence  of the mixed potential
function with the aid of which the system of differential
equations describing the dynamics of such a network is written
into a special form (see \S 4 in \cite{moser}). The mixed
potential function is constructed explicitly  only for the
networks whose graph possesses a tree containing all the capacitor
branches and none of the inductive branches, that is, the network
does no contain any loops of capacitors or cutsets of inductors,
each resistor tree branch corresponds to a current-controlled
resistor, each resistor co-tree branch corresponds to a
voltage-controlled resistor (see \S13 in \cite{moser}).

\noindent In \cite{maschke}, the dynamics of a nonlinear LC circuit is shown to be of Hamiltonian nature with respect to a certain  Poisson 
bracket which may be degenerate, that is, non-symplectic. In this formalism,  the constitutive relations of the inductors and 
capacitors are used to define the Hamiltonian function in terms of capacitor charges and inductor fluxes, while the topological 
constraints of the network graph and  Kirchhoff's laws define the Poisson bracket on the space of capacitor charge and inductor flux 
variables.

\noindent But for all those formulations,  a certain topological assumption on the electrical circuit appears to be crucial, that is, the 
circuit is supposed to  contain neither loops of capacitors
 nor cutsets of inductors.

\noindent In \cite{schaft}, \cite{maschke2} and \cite{bloch} the Poisson bracket is replaced by the more general notion of a Dirac structure 
on a vector space, leading to implicit Hamiltonian systems. The
Hamiltonian function is the total electromagnetic energy of the
circuit and the vectorial state space is defined by the inductors'
fluxes and capacitors' charges. The Dirac structure on the state
space is obtained from  Kirchhoff's laws. In this formalism, it is
possible to include networks which do not obey the topological
assumption mentioned before.

In the paper at hand we will see that the restricted class of
networks involving capacitor loops and inductor cut sets are
naturally captured  by the Birkhoffian approach.
We are going to discuss explicit examples  in order to demonstrate the advantages of the 
Birkhoffian approach in  the analysis of the resulting systems.
Another advantage of the Birkhoffian approach is the possible
inclusion of dissipative effects caused by resistors included in a
network.
It is a straight-forward matter to extend the 
approach presented here to the case of RLC circuits, that is circuits containing resistors in addition to capacitors and inductors. 
However, to start it appears to be more instructive to restrict
the theory to the case of LC circuits. The investigation of RLC
circuits will be presented into another paper.

The following parts of the paper are organized as follows. In Section 2 we recall the basics of Birkhoffian systems (see \cite{bir}) 
presented from the point view of differential geometry using the formalism of jets (see \cite{oliva}).
Birkhoffian formalism is  a global formalism of the dynamics of implicit systems of second order ordinary differential equations on 
a manifold.
In particular, we extend the approach in  \cite{oliva} to non-autonomous systems in order to be able to treat the case of networks 
with independent voltage and current sources later on in section
4. In Section 3, our Birkhoffian formulation of the dynamic
equations of a nonlinear LC circuit is introduced. Properties of
the corresponding Birkhoffian such as its regularity and its
conservativeness are also discussed in this section.  For a
nonlinear  LC electric network  each Birkhoffian of the family is
 conservative.
If there exists in the network some loop which contains only
capacitors the Birkhoffian is never regular. For such electrical
circuits, we present a systematic procedure to reduce the original
configuration space to a lower dimensional one, thereby
regularizing the Birkhoffian. On the reduced configuration space
the reduced Birkhoffian will still be conservative. In case the LC
circuit has loops which contain only linear inductors, the
original configuration space can be further reduced to a lower
dimensional one. Inductor loops can be regarded as some
conservative quantities of the network.
 In Section 4 we give a Birkhoffian formulation of a nonlinear
 LC circuit with independent sources and discuss in this context
 the concepts of regularity and conservativeness.
For instance, it turns out that voltage sources do not destroy
conservativeness, even in the nonlinear case, while current
sources might do so. Finally, in Section 5 we consider some
specific examples.
 These examples are supposed to serve our purpose of demonstrating
the power of the Birkhoffian approach. In particular, we can allow
 capacitor loops as well as  inductor cutsets, as already
mentioned before. Also, we investigate the question of
conservativeness of the underlying Birkhoffian in case of a
circuit with independent current and voltage sources.
\\

\textit{\textbf {Acknowledgement.}} We are grateful to Professor Marsden for a fruitful discussion concerning the topic of this 
paper.

\section{Birkhoffian systems}

For a smooth m-dimensional differentiable connected manifold $M$, we consider the tangent bundles $(TM,\pi_M,M)$ and ($TTM, 
\pi_{TM}, TM$).  Let $q=(q^1$, $q^2,$..., $q^m)$
be a local coordinate system on $M$.
This induces   natural local coordinate systems  on $TM$ and  $TTM$, denoted
by ($q,\dot{q}$), respectively ($q,\, \dot{q},\,dq,\, d\dot{q}$).
The \textit{2-jets manifold} $J^2(M)$ is a
 $3m$-dimensional submanifold of $TTM$ defined by
 \be
J^2(M)=\left\{ z\in TTM\ /\ T\pi_{M}(z)=\pi_{TM}(z)\right\} \label{0}
\ee
where $T\pi_{M}:TTM\to TM$ is the tangent map of $\pi_{M}$.
We write $\pi_J:=\pi_{TM}\arrowvert _{J^2(M)}=T\pi_M\arrowvert_{J^2(M)}$.
($J^2(M),\, \pi_J,\, TM$), called the \textit{2-jet bundle} (see \cite{oliva}), is an affine bundle modelled on the vertical vector 
bundle ($V(M),\, \pi_{TM}\arrowvert_{V(M)},\, TM$),  $
V(M)=\bigcup_{v\in TM} V_v(M) $ , where  $V_v(M)=\{z\in T_vTM\,
\arrowvert\, (T\pi_M)_v(z)=0\}$.
In \cite{marsden}, \cite{saunders} this bundle is denoted by $T^2(M)$ and 
named \textit{second-order tangent bundle}.
In natural local coordinates, the equality in (\ref{0}) yields
 $(q,\, \dot{q},\,\dot{q},\, d\dot{q}\arrowvert_{J^2(M)})$ as a local coordinate system on $J^2(M)$. We set $\ddot{q}:=d\dot 
{q}\arrowvert_{J^2(M)}$.
 Thus, a local coordinate system $q$ on $M$ induces the  natural local  coordinate system
 $(q,\, \dot{q},\,\ddot{q})$ on $J^2(M)$. For further details on this affine bundle see  \cite{marsden}, \cite{oliva}, 
\cite{saunders}.

A \textit{\textbf {Birkhoffian}} corresponding to the configuration manifold $M$ is a smooth 1-form $\omega$ on $J^2(M)$ such that, 
for any $x\in M$, we have
\be
i_{x}^*\omega=0\label{defbir}
\ee
where $i_x:\beta^{-1}(x)\to J^2(M)$ is the embedding of the submanifold $\beta^{-1}(x)$ into $J^2(M)$, $\beta=\pi_M\circ \pi_J$. 
From this definition it follows that, in the  natural local coordinate system ($q,\, \dot{q},\, \ddot{q}$) of $J^2(M)$, a 
Birkhoffian $\omega$ is given by \be \omega=\sum^{m}_{j=1}Q_j(q,\,
\dot{q},\, \ddot{q})dq^j\label{bircor} \ee with certain functions
$Q_j:J^2(M)\to \mathbf{R}$.

\noindent The pair ($M,\, \omega$) is said to be a \textit{
\textbf {Birkhoff system}} (see \cite{oliva}).

\noindent The \textit{\textbf{differential system associated to a
Birkhoffian}} $\omega $ (see \cite{oliva} ) is  the set (maybe
empty)  $D(\omega$), given by \be D(\omega):=\left\{z\in J^2(M)\
\arrowvert\, \omega(z)=0 \right\}\ee The manifold $M$ is the
\textit{space of configurations} of $D(\omega)$, and $D(\omega)$
is said to have $m$ 'degrees of freedom'. The $Q_i$ are the
'generalized external forces' associated to the local coordinate
system $(q)$. In the natural local coordinate system, $D(\omega)$
is characterized by the following implicit system of second order
ODE's \be Q_j(q,\, \dot{q},\, \ddot{q})=0 \textrm{ for all }
j=\overline{1,m}\label{difsistem} \ee \noindent {\it We conclude
that the Birkhoffian formalism is  a global formalism for the
dynamics of implicit systems of second order  differential
equations on a manifold}.
\\

 Let us now associate a vector field to a Birkhoffian $\omega$.\\
 A \textit{vector field} $Y$ on the manifold $TM$ is a smooth
function $Y:TM\to TTM$ such that $\pi_{TM}\circ Y$=id.
Any vector field $Y$ on $TM$ is called a \textit{second order vector field on TM} if and only if $T{\pi_M}(Y_v)=v$ for all $v\in 
TM$.

\noindent A cross section $X$ of the affine bundle ($J^2(M),\, \pi_J,\, TM$), that is, a smooth function $X:TM\to J^2(M)$ such that  
$\pi_J\circ X$=id,
 can be identified with a special vector field on $TM$, namely, the second order vector field on $TM$
 associated to $X$.  Indeed, because ($J^2(M),\, 
\pi_J,\, TM$) is a sub-bundle of
($TTM, \pi_{TM}, TM$) as well as of ($TTM, T\pi_{M}, TM$), its sections can be regarded as sections of these two tangent bundles. Thus, using 
the canonical   embedding  $i:J^2(M)\to TTM$,  $X$ can be identified with $Y$, that is,  $Y=i\circ X$.

\noindent In  natural  local coordinates a second order vector
field can be represented as \be
Y=\sum^{m}_{j=1}\left[\dot{q}^i\f{\pa}{\pa q^i}+\ddot{q}^i(q,\,
\dot{q})\f{\pa}{\pa \dot{q}^i}\right]\label{vectorfield} \ee

A \textit{\textbf {Birkhoffian vector field}} associated to a Birkhoffian $\omega$ of $M$ (see \cite{oliva}) is a smooth second 
order vector field on $TM$, $Y=i\circ X$, with $X:TM\to J^2(M)$,
such that $Im X\subset\, D(\omega)$, that is \be X^*\omega=0 \ee

\noindent In the natural  local coordinate system, a Birkhoffian vector field is given by  (\ref{vectorfield}), such 
that $Q_j(q,\dot{q},\ddot{q}(q,\dot{q}))=0$.

A Birkhoffian $\omega $ is \textit{\textbf {regular}} if and only
if  \be \textrm{det}\left[\f{\pa Q_j}{\pa \ddot{q}^i}(q,\,
\dot{q},\, \ddot{q})\right]_{i,j=1,...,m}\neq 0\label{regular} \ee
for all $(q,\, \dot{q},\, \ddot{q})$, and for each $(q,\, \dot{q})$, there exists $\ddot{q}$ such that $Q_j(q,\, \dot{q},\, 
\ddot{q})=0,\, j=1,...,m.$

\noindent If a Birkhoffian $\omega$ of $M$ is regular, then it satisfies \textit{the principle of determinism}, that is, there exists an 
unique Birkhoffian vector field $Y=i\circ X$ associated to $\omega$ such that $Im\, X=D(\omega)$ (see \cite{oliva}).

 A Birkhoffian $\omega$ of $M$ is called \textit{\textbf{conservative}} if and only if there exists a smooth function 
$E_{\omega}:TM\to \mathbf{R}$ such that

\be
 (X^*\omega)Y=dE_{\omega}(Y)\label{conserv'}
\ee
for all second order vector fields $Y=i\circ X$ (see \cite{oliva}).

\noindent Equation (\ref{conserv'}) is  equivalent, in the natural
local coordinate system, to the identity (see \cite{bir}, p. 16,
eq. 4) \be
\sum^{m}_{j=1}Q_j(q\, \, \dot{q},\, \ddot{q})\dot{q}^j=\sum^{m}_{j=1}\left[\f{\pa E_{\omega}}{\pa q^j}\dot{q}^j+\f{\pa E_{\omega}}{\pa 
\dot{q}^j}\ddot{q}^j\right]\label{conserv} \ee

 \noindent $E_{\omega}$ is constant on $TM$ if and only if $dE_{\omega}(Y)=0$
for all second order vector fields $Y$ on $TM$ (see \cite{oliva}).

\noindent If $\omega$ is conservative and $Y$ is a Birkhoffian
vector field, then (\ref{conserv'}) becomes \be dE_{\omega}(Y)=0
\ee This means that $E_{\omega}$ is constant along the
trajectories of $Y$.

It is also possible to introduce, in a natural manner, the notion
of constrained Birkhoff system (see \cite{oliva}, \S 4).

\noindent Let  ($M$, $\omega$) be  a Birkhoff system and
$\mathfrak{S}$ a smooth constant rank  affine sub-bundle of the
affine bundle $\pi_{J}:J^2(M)\longrightarrow TM.$
 Locally, the submanifold $\mathfrak{S}$ of codimension
$\textsc{n}$, is described by the vanishing of $\textsc{n}$
independent affine functions \be
\phi^{\nu}(q,\dot{q},\ddot{q})=\sum^{m}_{i=1}b^\nu_i(q,\dot{q})\ddot{q}^i+a^\nu(q,\dot{q}),
\quad \nu=\overline{1,\textsc{n}} \ee

\noindent A triple ($M$, $\omega$, $\mathfrak{S}$) is called
\textbf{\textit{constrained Birkhoff system}}.

\noindent The \textbf{\textit{constrained differential system
associated to }} the constrained Birkhoff system ($M$, $\omega$,
$\mathfrak{S}$)  is the set \be D(\omega, \mathfrak{S})= \{z\in
\mathfrak{S}\arrowvert\,\omega(z)=0 \} \label{c1} \ee
\\

 Let us now generalize these concepts to time-dependent dynamical
systems.

\noindent For the usual formulation of Lagrangian and Hamiltonian time-dependent mechanics (see for example \cite{marsden}, \S 5.1, 
\cite{saunders} \S 4.1, \S 4.6
), the configuration space has the form $\mathbf{R}\times M$, the phase space has the form $\mathbf{R}\times T^*M $, and the 
velocity space has the form $\mathbf{R}\times TM$, with some manifold M. If $(t,q)$ is a coordinate system on $\mathbf{R}\times M$, 
then  $(t,q, \dot{q})$ is a  coordinate  system on $\mathbf{R}\times TM$.
Thus, $\mathbf{R}\times TM$ can be interpreted  as a submanifold of $T(\mathbf{R}\times M)$ given by
\be
\dot{t}=1\label{timp1}
\ee
 From the physical point of view,  this means that a reference frame has been chosen. This is not the case for relativistic 
mechanics. The reference system provides a splitting between the
time and the state coordinates of a mechanical system. Within the
Birkhoffian framework, we follow  the usual non-relativistic
lines. Thus,
 for the time-dependent system, we have in addition the equation
\be
\ddot{t}=0\label{timp2}
\ee

\noindent In view of (\ref{timp1}), (\ref{timp2}), we choose in
the study of time-dependent dynamical systems the extended bundle
$\mathbf{R}\times J^2( M)$.

A {\bf \textit{ time-dependent Birkhoffian}} is a smooth family of
1-forms $\omega_{t}$ on $J^2( M)$ defined by \be
\omega_{t}=\sum^{m}_{j=1}Q_j(t,q, \dot{q},
\ddot{q})dq^j\label{bircortime} \ee where   ($t,q,\dot{q},
\ddot{q}$) is the natural  local coordinate system on
$\mathbf{R}\times J^2( M)$. Thus, our time-dependent Birkhoffian
is obtained by merely freezing $t$ and constructing the
Birkhoffian for any fixed
value of  $t$ as before.\\
A \textit{time-dependent second order vector field } (see
\cite{saunders}) on $\mathbf{R}\times TM$ has the following
representation in the  natural  local coordinate system \be
Y_{t}=\f {\pa}{\pa t}+\sum^{m}_{j=1}\left[\dot{q}^j\f{\pa}{\pa
q^j}+\ddot{q}^j(t,q, \dot{q})\f{\pa}{\pa
\dot{q}^j}\right]\label{yt} \ee

\noindent Thus, for a time-dependent system, a {\bf \textit {time-dependent  Birkhoffian  vector field}} on $\mathbf{R}\times T( M)$ has the 
representation (\ref{yt}), where
 $Q_j(t,q,\dot{q},\ddot{q}(t,q,\dot{q}))=0$.

 A time-dependent Birkhoffian $\omega_{t}$ is \textit{\textbf
{regular}} if and only if  \be \textrm{det}\left[\f{\pa Q_j}{\pa
\ddot{q}^i}(t,q,\, \dot{q},\, \ddot{q})\right]_{i,j=1,...,m}\neq
0\label{regularna} \ee
for all $(t,q,\, \dot{q},\, \ddot{q})$, and for each $(t,q,\, \dot{q})$, there exists $\ddot{q}$ such that $Q_j(t,q,\, \dot{q},\, 
\ddot{q})=0,\, j=1,...,m.$

 A time-dependent Birkhoffian  $\omega_{t}$ is called
\textit{\textbf{conservative}} if  and only if
there exists a smooth family of functions $E_{\omega_{t}}: 
TM \longrightarrow \mathbf{R}$ such that, everywhere, \ba
\sum^{m}_{j=1}Q_j(t,q\, \, \dot{q},\, \ddot{q})\dot{q}^j=
\sum^{m}_{j=1}\left[\f{\pa E_{\omega_{t}}}{\pa q^j}\dot{q}^j+
\f{\pa E_{\omega_{t}}}{\pa \dot{q}^j}\ddot{q}^j\right]
\label{surseconserv} \ea

\noindent If $\omega_{t}$ is conservative and $Y_{t}$ is a
time-dependent Birkhoffian vector field then, from
(\ref{surseconserv}), we obtain the \textbf{\textit{generalized
balance law}} \be \f{dE_{\omega_{t}}}{dt}=\f {\pa
E_{\omega_{t}}}{\pa t} \ee along trajectories of $Y_{t}$.

\section{LC circuit dynamics}

A simple electrical circuit provides us with an \textit{oriented
connected} graph, that is, a collection of points, called nodes,
and a set of connecting lines or arcs, called branches, such that
in each branch is given a direction and there is at least one path
 between any two
nodes. A path is a sequence of branches such that the origin of
the next branch coincides with the end of the previous one. The
graph will be assumed to be \textit{planar}, that is, it can be
drawn in a plane without branches crossing. For the graph theoretic terminology, see, for 
example \cite{graph}.\\
 Let $b$  be the total number of branches in the graph,
$n$ be one less than the number of  nodes and $m$
be the cardinality of a 
selection of loops that cover the whole graph. Here, a loop is a
path such that the first and last node coincide and that does not
use the same branch  more than once. By Euler's polyhedron
formula, $b=m+n$.\\
A \textit{cutset} in a connected graph, is a minimal set of branches whose  removal from the graph, renders the graph 
disconnected. For example, the set branches tied to a node is a
cutset.

\noindent
 We choose a reference node and a  current direction in
each   $l$-branch of the graph,
 $l=1,...,b$. We also consider a covering of the graph with $m$ loops,
   and a current direction in  each $j$-loop, 
$j=1,...,m$. We assume that the associated graph has at least one
loop, meaning that $m>0$.\\
 A graph  can be described by matrices: a ($bn$)-matrix
$B\in \mathfrak{M}_{bn}(\mathbf{R})$, rank$(B)=n$, called
\textit{incidence matrix} and a ($bm$)-matrix $A\in \mathfrak{M}_
{bm} (\mathbf{R})$, rank$(A)=m$, called \textit{loop matrix}.
These matrices contain only 0,1, -1. An element of the matrix $B$
is 0 if a branch  $b$ is not incident with a node $n$, 1 if branch
$b$ enters  node  $n$ and -1 if
 branch $b$ leaves node $n$, respectively. An element of the matrix $A$ is 0 if a branch
 $b$ does
not belong to a loop $m$, 1 if  branch $b$
 belongs to loop $m$
and their directions agree
 and -1 if
 branch $b$ belongs to loop $m$  and their  directions
oppose, respectively. For the fundamentals of electrical circuit
theory, see, for example \cite{chua}.

 The states of the circuit have two components, the
currents through the branches, denoted by $\textsc{i}\in
\mathbf{R}^b$,  and the voltages across the branches, denoted by
$v\in \mathbf{R}^b$.
 Using the matrices
$A$ and $B$, Kirchhoff's current law and Kirchhoff's voltage law
can be expressed by the equations \be B^T \textsc{i} =0 \quad
(KCL)\label{5} \ee \be A^T v =0\quad (KVL)\label{6} \ee \noindent
Tellegen's theorem establishes a relation between the matrices
$A^T$ and $B^T$: \textit{the kernel of the matrix $B^T$ is
orthogonal to the kernel of the matrix $A^T$}(see e.g.,
\cite{moser} page 5).

 The next step is to introduce the branch elements in a
simple electrical circuit. The branches of the
 graph associated to an  LC electrical circuit, can be classified
into two categories:  inductor branches
 and capacitor branches. A capacitor loop will contain only
 capacitor branches and an inductor cutset will  contain only inductor branches. Let
 $k$ denote the number of inductor branches and $p$ the number of capacitor
 branches, respectively. We assume that just one electrical device is associated to each
branch, then, we have $b=k+p$. Thus, we can write
$(\textsc{i}_{a},\textsc{i}_{\alpha})\in \mathbf{R}^r\times
 \mathbf{R}^p\simeq\mathbf{R}^b$, where
 $\textsc{i}_ {a}$, $\textsc{i}_{\alpha}$ are the currents
through the inductors,  the capacitors, respectively, and
$v=(v_{a},v_{\alpha})\in
\mathbf{R}^k\times\mathbf{R}^p\simeq\mathbf{R}^b$, where $v_{a}$,
$v_\alpha$ describe the voltage drops across the
 the inductors, the capacitors, respectively.\\
To exemplify, let us now write the matrices $B$ and $A$, for a
circuit which contains four inductors, three capacitors and which
has the following oriented connected graph
\begin{center}
\scalebox{0.40} {\includegraphics*{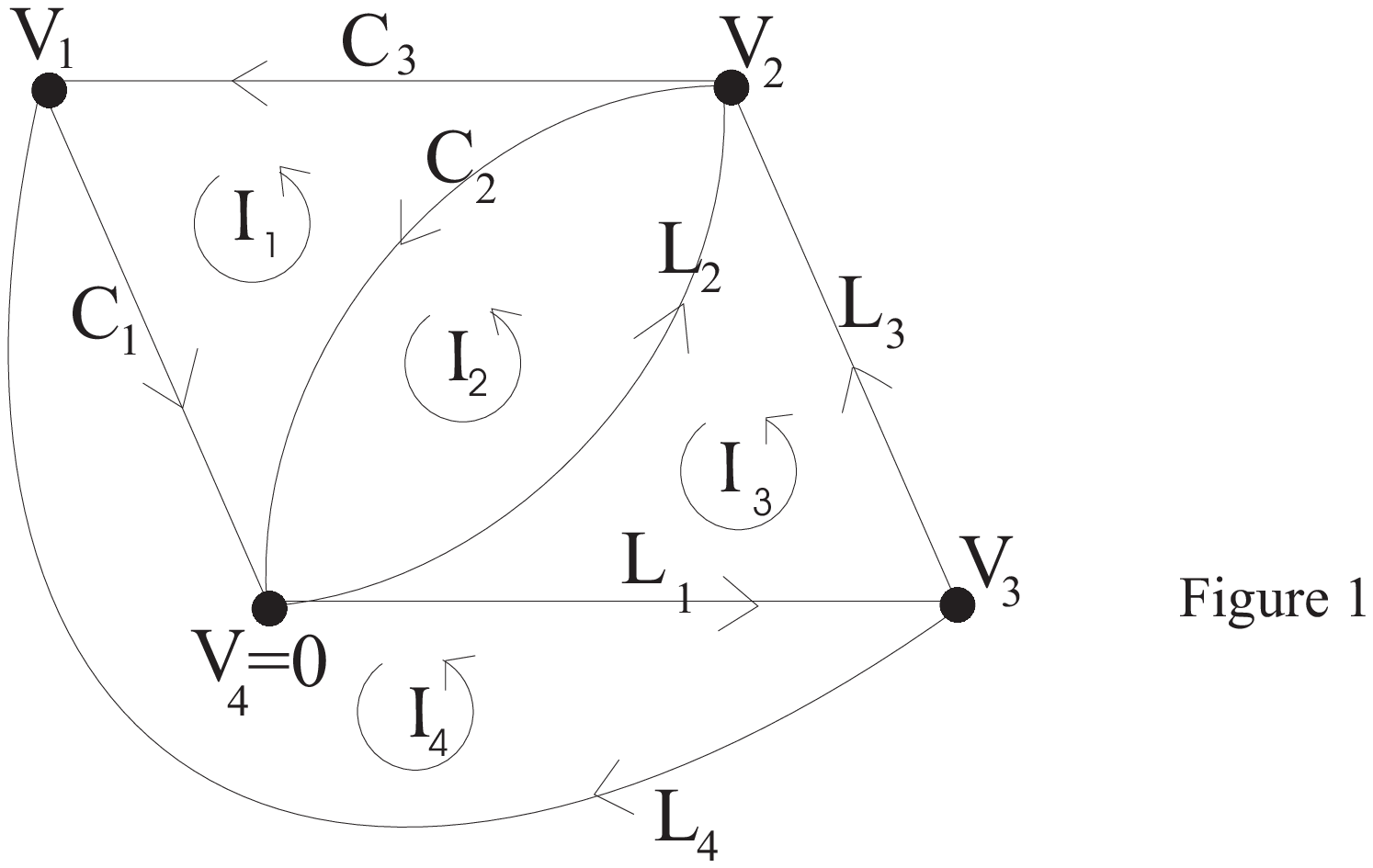}}
 \label{figura}
\end{center}
\noindent We have $k=4$, $p=3$, $n=3$, $m=4$ , $b=7$.
 We
choose the reference node to be $V_4$ and the current directions
as indicated in Figure 1. We cover the  graph with the loops $I_1,
\, I_2,\, I_3,\, I_4$. Let $V=(V_1,V_2,V_3)\in \mathbf{R}^3$ be
the vector of node voltage values,
 $\textsc{i}=(\textsc{i}_{a},\textsc{i}_\alpha)\in
\mathbf{R}^4\times \mathbf{R}^3$ be
 the vector of branch current  values  and
 $v=(v_{a},v_{\alpha})\in \mathbf{R}^4\times\mathbf{R}^3$ be
 the vector of branch  voltage
 values.\\
The branches in Figure 1 are  labelled as follows: the first, the
second, the third and the fourth branch are the inductor branches
$\textsc{L}_1$, $\textsc{L}_2$, $\textsc{L}_3$, $\textsc{L}_4$ and
the last three branches are the capacitor branches $\textsc{C}_1$,
$\textsc{C}_2$, $\textsc{C}_3$. The incidence and loop matrices,
$B\in \mathfrak{M}_{73}(\mathbf{R})$ and
 $A\in \mathfrak{M}_{74}(\mathbf{R})$, write as
 \be
 B=\left(
\begin{array}{ccc}
 0& 0&1\\
 0&1& 0\\
 0&1& -1\\
1& 0& -1\\
 -1& 0& 0 \\
 0& -1& 0 \\
1& -1& 0
\end{array}
\right), \quad
 A=\left(
\begin{array}{cccc}
 0& 0 & 1 & -1 \\
 0& 1 & -1 & 0 \\
 0& 0 & 1 & 0\\
 0& 0 & 0 & -1\\
 1& 0 & 0 & -1\\
-1& 1 & 0 & 0 \\
 1& 0&  0& 0
\end{array}
\right)
 \label{ex91}
\ee For another choice of the covering loops and of the current
directions in the loops we obtain a different matrix $A$
and for another choice of the reference node and of the current
directions in the branches we  obtain a different matrix $B$.
$\quad\blacksquare$
\\

\noindent Each capacitor is supposed to be charge-controlled. For
the nonlinear capacitors we assume \be
v_\alpha=C_\alpha(\textsc{q}_{\alpha}), \quad \alpha=1,..., p
\label{2} \ee where the functions
$C_\alpha:\mathbf{R}\longrightarrow \mathbf{R}\backslash \{0\}$
are smooth and invertible, and the 
$\textsc{q}_\alpha$'s denote the charges of the capacitors.
 The current through a  capacitor is given by  the
time-derivative
 of the corresponding  charge
\be \textsc{i}_\alpha=\f{d\textsc{q}_{\alpha}}{dt}, \quad
\alpha=1,...,p \label{4'} \ee $t$ being the time variable.\\
 Each inductor is supposed to be current-controlled. For
the  nonlinear inductors we assume \be
v_a=L_{a}({\textsc{i}_a})\f{d{\textsc{i}_a}}{dt}, \quad a=1,...,k
\label{4} \ee where $L_a:\mathbf{R}\longrightarrow
\mathbf{R}\backslash \{0\}$ are smooth invertible functions.\\

\noindent If the capacitors and the inductors are linear then the
relations (\ref{2}) and (\ref{4}) become, respectively, \be
v_\alpha=\frac{\textsc{q}_\alpha}{\textrm {\scriptsize C}_\alpha},
\quad v_a=\textrm {\scriptsize
L}_a\f{d{\textsc{i}_a}}{dt}\label{15} \ee where {\scriptsize
C}$_\alpha\neq 0$ and {\scriptsize L}$_a\neq 0$ are distinct
constants.
\\

 Taking into account (\ref{2}), (\ref{4'}), (\ref{4}),
the equations (\ref{5}), (\ref{6}), become
 \be \left\{
\begin{array}{ll}
B^T
\left(
\begin{array}{c}
\textsc{i}_a \\
\f{d{\textsc{q}}_{\alpha}}{dt}
\end{array}
\right)=0\\
\\
A^T
\left(
\begin{array}{c}
L_{a}(\textsc{i}_a)\,
\f{d{\textsc{i}_a}}{dt} \\
C_\alpha(\textsc{q}_{\alpha})
\end{array}
\right)=0
\end{array}
\right.\label{7} \ee

\noindent  \textit{In the following we give a Birkhoffian
formulation for the network described by the system of equations
(\ref{7}). Using the first set of equations (\ref{7}),
 we are going to define a  family of
  $m$-dimensional affine-linear configuration spaces
  $M_c\subset\mathbf{R}^b$ parameterized by a constant vector
  $c$ in $\mathbf{R}^n$. This vector is related to the initial
  values of the $\textsc{q}$-variables at some instant of time.
  At this point we notice that actually already the initial values 
corresponding to the $\textsc{q}$-variables associated to capacitors,
together with those of $m$ distinguished branch currents 
denoted by $\dot{q}^j$ below, parameterize the whole solution set
of the equations in (\ref{7}). A Birkhoffian $\omega_c$ of the
configuration space $M_c$ arises from a linear combination of the
second set of equations (\ref{7}). Thus, ($M_c,\omega_c$) will be
a family of Birkhoff systems that describe the  LC circuit
considered.}
\\

 \noindent We  notice that the first set of equations (\ref{7})
 remains exactly the same for linear and nonlinear electrical devices.
 Thus, for obtaining the  configuration space, it is not important
 whether the devices are linear or nonlinear.
We shall see below
that  the only difference is that one ends up with a nonlinear configuration space
or rather configuration manifold when one 
regularizes the resulting Birkhoffian system in the case of
nonlinear networks.
\\

 Let  $H: \mathbf{R}^b\longrightarrow\mathbf{R}^n$ be a
linear map that, with respect to a  coordinate system
($x^1,...,x^b$) on $\mathbf{R}^b$, is given  by \be
H(x^1,...,x^b)=B^T \left(
\begin{array}{c}
x^1 \\
\vdots\\
x^b
\end{array}
\right)\label{9} \ee Then, $H^{-1}(c)$, with $c$ a constant vector
in $\mathbf{R}^n$, is an affine-linear subspace in $\mathbf{R}^b$.
Its dimension is $m=b-n$, because rank$(B)=n$.
\\

\noindent We define  $M_c$ as \be M_c:=H^{-1}(c)\label{8} \ee We
denote a  coordinate system on $M_c$ by  $q=(q^1,...,q^m)$. Then,
the natural coordinate system on the 2-jet bundle $J^2(M_c)$ is
$(q,\dot{q},\ddot{q})$.

\noindent Let us now represent the Birkhoffian in a specific
coordinate system on $M_c$:\\
 In the vector space $\mathbf{R}^k$,
 we identify points and vectors \be
\textsc{i}_a:=\f{d{\textsc{q}}_{(a)}}{dt},\label{var} \ee where
$(\textsc{q}_{(a)})_{a=1,...,k}$ is a  coordinate system on
$\mathbf{R}^k$. Taking into account (\ref{var}) and the fact that
the matrix $B^T$ is a constant matrix, we integrate  the first set
of equations (\ref{7}) to arrive at \be B^T \left(
\begin{array}{c}
\textsc{q}_{(a)} \\
\textsc{q}_{\alpha}
\end{array}
\right)=c\label{7''} \ee with $c$ a constant vector in
$\mathbf{R}^n$.\\
 Likewise consider coordinates on
$\mathbf{R}^b\simeq \mathbf{R}^k\times\mathbf{R}^p$ defined by \be
x^1:=\textsc{q}_{(1)},\,
...,\,x^k:=\textsc{q}_{(k)},\,x^{k+1}:=\textsc{q}_{1},\,...,\,
x^b:=\textsc{q}_{p}\label{cs} \ee From  (\ref{9}), (\ref{8}), we
see that we can define coordinates on $M_c$ by solving the
equations in
(\ref{7''}) in terms of an 
appropriate set of  $m$ of the \textsc{q}-variables, say $q=(q^1,...,q^m)$.
In other words, we express any of the x-variables  as a function of $q=(q^1,...,q^m)$,
 namely,
\ba &&x^a= \sum_{j=1}^{m}\mathcal{N}^{a}_{ j}q^j+const,\,
a=1,...,k,\nonumber\\
&&x^{\alpha}= \sum_{j=1}^{m}\mathcal{N}^{\alpha}_{j}q^j+const, \,
\alpha=k+1,...,b \label{10'} \ea with certain constants
$\mathcal{N}^{a}_{j}$, and $\mathcal{N}^{\alpha}_{j}$. Here we can
think of the constants $const$ as being initial values of the
$x$-variables at
 some instant of time.

\noindent From  (\ref{4'}), (\ref{var}),
 (\ref{cs}) and  differentiating  (\ref{10'}) we get
\be \textsc{i}=\mathcal{N}\dot{q}\label{N} \ee with the matrix of
constants $\mathcal{N}\in\mathfrak{M}_{bm} (\mathbf{R})$, for some
$\dot{q}\in \mathbf{R}^m$.

\noindent Using  Tellegen's theorem and a fundamental theorem of
linear algebra, we now find a relation between the matrices
$\mathcal{N}$ and $A$. By a fundamental theorem of linear algebra
we have \be (Ker(A^T))^\perp=Im(A) \label{linearalgebra}\ee where
$A\in \mathfrak{M}_{bm}(\mathbf{R})$, $Ker(A^T):=\{x\in
\mathbf{R}^b\, |\, A^Tx=0\}$ is the kernel of $A^T$,
$Im(A):=\{x\in \mathbf{R}^b\, |\, Ay=x,\, \textrm{for some} \,
y\in \mathbf{R}^m \}$ is the image of $A$ and  $^\perp$ denotes
the orthogonal complement in $\mathbf{R}^b$ of
the respective vector subspace.\\
 For the incidence  matrix $B\in
\mathfrak{M}_{bn}(\mathbf{R})$ and the loop matrix $A\in
\mathfrak{M}_{bm} (\mathbf{R})$, which satisfy Kirchhoff's law
(\ref{5}), (\ref{6}), Tellegen's theorem writes as \be
Ker(B^T)=(Ker(A^T))^\perp \label{tellegen}\ee From the first set
of equations in (\ref{7}), and by constraction of the matrix
$\mathcal{N}$ in (\ref{N}), we have \be
Ker(B^T)=Im(\mathcal{N})\label{N1}\ee Therefore, using
(\ref{linearalgebra}), (\ref{tellegen}), (\ref{N1}), we obtain
$Im(A)=Im(\mathcal{N})$. Then, another application of
(\ref{linearalgebra}) yields \textit{ \be
Ker(A^T)=Ker(\mathcal{N}^T)\label{ker}\ee} Taking into account
(\ref{ker}), we see that there exists a nonsingular matrix
$\mathcal{C}\in \mathfrak{M}_{mm}(\mathbf{R})$ satisfing \be
\mathcal{C}A^T=\mathcal{N}^T\label{equivalence} \ee The matrix
$\mathcal{C}$ provides a relation between the  vector of the $m$
independent loop currents and the coordinate vector $q$ introduced
on  $M_c$.
\\

 \textit{Taking into account (\ref{ker}), we define the
Birkhoffian $\omega_c$ of $M_c $ such that the differential system
(\ref{difsistem}) is the linear combination of the second set of
equations in (\ref{7}) obtained by replacing $A^T$ with the matrix
$\mathcal{N}^T$. Thus, in terms of q-coordinates as chosen before,
the expressions of the components
$Q_j(q, \, \dot{q},\,\ddot{q})$ from  
(\ref{bircor}) are \be Q_j(q, \,
\dot{q},\,\ddot{q})=F_j(\dot{q})\ddot{q}+G_j(q), \quad
j=1,...,m\label{bir} \ee where \be
F_j(\dot{q})\ddot{q}=\sum^{k}_{a=1}\mathcal{N}^{\,  a}_{j}
L_a\left( \sum^{m}_{l=1}\mathcal{N}^ {a}_{\,
l}\dot{q}^l\right)\sum ^{m}_{i=1}\mathcal{N}^ {a}_{\,
i}\ddot{q}^i=\sum^{m}_{i=1} \left(\sum^{k}_{a=1}\mathcal{N}_j^{\,
a}\mathcal{N}^ {a}_{\,
i}\widetilde{L}_a\left(\dot{q}\right)\right)\ddot{q}^i
\label{bir2} \ee \be G_j(q)=\sum
^{b}_{\alpha=k+1}\mathcal{N}_j^{\,
\alpha}C_{\alpha-k}\left(\sum^{m}_{l=1}\mathcal{N}^ {\alpha}_{\,
l}q^l+const\right)= \sum ^{b}_{\alpha=k+1}\mathcal{N}_j^{\,
\alpha}\widetilde{C}_{\alpha-k}\left(q\right)\label{bir1} \ee}

 \textit{We claim that the Birkhoffian (\ref{bir})  is a {\bf conservative}
one}.
\\

\noindent Indeed, for our problem, the relation (\ref{conserv})
becomes \be
\sum^{m}_{j=1}\left[\left(\sum^{m}_{i=1}\sum^{k}_{a=1}\mathcal{N}^{\,
a}_{j}\mathcal{N}^{a}_{\,  i}\widetilde{L}_a(\dot{q})
\ddot{q}^i\right)\dot{q}^j+G_j(q)\dot{q}^j\right]=\sum^{m}_{j=1}\left[
\f{\pa E_{\omega_c}}{\pa q^j}\dot{q}^j+\f{\pa E_{\omega_c}}{\pa
\dot{q}^j}\ddot{q}^j\right] \label{bir3} \ee or (changing the
indices of summation) \be
\sum^{m}_{j=1}\left[\left(\sum^{m}_{i=1}\sum^{k}_{a=1}\mathcal{N}^{\,
a}_{i}\mathcal{N}^ {a}_{\,  j}\widetilde{L}_a\left(\dot{q}\right)
\dot{q}^i\right)\ddot{q}^j+G_j(q)\dot{q}^j\right]=\sum^{m}_{j=1}\left[\f{\pa
E_{\omega_c}}{\pa q^j}\dot{q}^j+ \f{\pa E_{\omega_c}}{\pa
\dot{q}^j}\ddot{q}^j\right] \label{bir4} \ee

\noindent Because of the special form of the terms on the left hand side of (\ref{bir4}), we can look for the required function 
$E_{\omega_c}(q,\dot{q})$ as a sum of a function  depending only
on $q$, and a function  depending only
on $\dot{q}$. From the theory of 
total differentials, a necessary condition for the existence of
such functions is the fulfilment of the following relations \be
\left\{
\begin{array}{ll}
\f {\partial G_j(q)}{\partial q^l}-\f {\partial G_l(q)}{\partial q^j}=0\\
\\
\f {\partial \mathcal{F}_j(\dot{q})}{\partial \dot{q}^l }-
\f {\partial \mathcal{F}_l (\dot{q})}{\partial \dot{q}^j}=0
\end{array}
\right.\label{bir5} \ee for any $j,l=1,...,m$, where \be
\mathcal{F}_j(\dot{q}):=\sum^{m}_{i=1}\sum^{k}_{a=1}\mathcal{N}^{\,
a}_{i}\mathcal{N}^ {a}_{j}\widetilde{L}_a\left(\dot{q}\right)
\dot{q}^i \label{bir6} \ee In view of (\ref{bir1}), (\ref{bir6})
we get: \be \f {\partial G_j(q)}{\partial q^l}= \sum
^{b}_{\alpha=k+1}\mathcal{N}^{\alpha}_j\mathcal{N}^{\alpha}_{l}\widetilde{C}_{\alpha-k}'(q)
\label{bir7} \ee \be \f {\partial \mathcal{F}_j(\dot{q})}{\partial
\dot{q}^l}=
\sum^{k}_{a=1}\mathcal{N}^a_{l}\mathcal{N}^{a}_{j}\widetilde{L}_a(\dot{q})+
\sum^m_{i=1}\sum^{k}_{a=1}\mathcal{N}^a_{i}\mathcal{N}^
{a}_{j}\mathcal{N}^ {a}_{l}\widetilde{L}_a'(\dot{q})\dot{q}^i
\label{bir8} \ee
where $\widetilde{C}'_{\alpha}:=\f {d\widetilde{C}_{\alpha}(\eta)}{d\eta}$, $\widetilde{L}'_{a}:=\f 
{d\widetilde{L}_a(\eta)}{d\eta}$. Therefore, the left hand side of
  (\ref{bir5}) become \be \sum^{b}_{\alpha=k+1}
(\mathcal{N}^{\alpha}_j\mathcal{N}^{\alpha}_l-
\mathcal{N}^{\alpha}_l\mathcal{N}^{\alpha}_j)\widetilde{C}'_{\alpha-k}(q)\label{bir9}
\ee \be \sum^{k}_{a=1} \left(
\mathcal{N}^{a}_l\mathcal{N}^a_{j}-\mathcal{N}^{a}_j\mathcal{N}^{a}_l\right)\left(\widetilde{L}_{a}(\dot{q})
-\widetilde{L}'_{a}(\dot{q})(\sum^{m}_{i=1}\mathcal{N}_{i}^{a}\dot{q}^{i})\right)
\label{bir10} \ee We now easily see that the expressions in
(\ref{bir9}), (\ref{bir10}) are zero and (\ref{bir5}) are
satisfied. Thus, we proved the existence of a function
$E_{\omega_c}(q,\dot{q})$
 such that (\ref{bir4}) is fulfilled. \\
  Let us now look for the expression of this
 function.
For linear devices, taking into account (\ref{15}), we have \ba
\widetilde{L}_a(\dot{q})={\textrm {\scriptsize L}}_a, \quad
\widetilde{C}_{\alpha-k}(q)=\frac{\sum^{m}_{i=1}\mathcal{N}^
{\alpha}_{\,  i}q^i}{{\textrm {\scriptsize C}}_{\alpha-k}}+const
\ea with ${\textrm {\scriptsize L}}_a$, ${\textrm {\scriptsize
C}}_\alpha$ being real constants. Therefore, the functions
$\mathcal{F}_j(\dot{q})$ and $G_j(q)$ from (\ref{bir6}),
(\ref{bir1})  become \be
\mathcal{F}_j(\dot{q}):=\sum^{m}_{i=1}\sum^{k}_{a=1} {\textrm
{\scriptsize L}}_a \mathcal{N}^{\,  a}_{i}\mathcal{N}^ {a}_{j}
\dot{q}^i \ee \be G_j(q):=\sum ^{b}_{\alpha=k+1}
\sum^{m}_{i=1}\f{\mathcal{N}_j^{\,  \alpha}\mathcal{N}^
{\alpha}_{\,  i}}{{\textrm {\scriptsize
C}}_{\alpha-k}}q^i+(const)_j \ee Thus, in the linear case, it is
not difficult to find  the
 function $E_{\omega_c}(q,\dot{q})$ such that (\ref{bir4})  is satisfied. This is  \be
E_{\omega_c}(q,\dot{q})=\f{1}{2}\sum^{k}_{a=1}\sum_{i,j=1}^m
{\textrm {\scriptsize L}}_a
\mathcal{N}^{a}_i\mathcal{N}^{a}_j\dot{q}^i \dot{q}^j+\f{1}{2}
\sum^{b}_{\alpha=k+1}\sum^{m}_{i,j=1}
\f{\mathcal{N}^{\alpha}_i\mathcal{N}^{\alpha}_j} {{\textrm
{\scriptsize C}}_{\alpha-k}}q^iq^j+\sum_{j=1}^{m}(const)_jq^j
\label{energie} \ee In order to derive such a function for
nonlinear devices, we start with the  equations \be \left\{
\begin{array}{ll}
\f{\pa E_{\omega_c}}{\pa \dot{q}^1}=\mathcal{F}_1(\dot{q})=
\sum^{m}_{i=1}\sum^{k}_{a=1}\mathcal{N}^{\,  a}_{i}\mathcal{N}^
{a}_{1}\widetilde{L}_a\left(\dot{q}\right)
\dot{q}^i\\
\\
\f{\pa E_{\omega_c}}{\pa q^1}=G_1(q)=\sum
^{b}_{\alpha=k+1}\mathcal{N}_1^{\,
\alpha}\widetilde{C}_{\alpha-k}\left(q\right)
\end{array}
\right. \ee Integrating with respect to $q^1$ and $\dot{q}^1$,
respectively, we get \ba\hspace{-0.7cm}
E_{\omega_c}(q^1,...,q^m,\dot{q}^1,...,\dot{q}^m)&=&\sum_{a=1}^{k}\int
\widetilde{L}_a\left(\dot{q}\right)\mathcal{N}^{\,  a}_{i}
\dot{q}^i\mathcal{N}^
{a}_{1}d\dot{q}^1+f_1(\dot{q}^2,...,\dot{q}^m)+\nonumber\\
&&
\sum_{\alpha=k+1}^{b}\int \widetilde{C}_{\alpha-k}\left(q\right)\mathcal{N}_1^{\,  \alpha}dq^1+g_1(q^2,...,q^m)
\label{j=1}
\ea
$f_1$ depends only on $\dot{q}^2,...,\dot{q}^m$ and $g_1$ depends only on $q^2,...,q^m$.
For $j=2$, we have
\be
\left\{
\begin{array}{ll}
\f{\pa E_{\omega}}{\pa \dot{q}^2}=\mathcal{F}_2(\dot{q})=
\sum^{m}_{i=1}\sum^{k}_{a=1}\mathcal{N}^{\,  a}_{i}\mathcal{N}^
{a}_{2}\widetilde{L}_a\left(\dot{q}\right)
\dot{q}^i\\
\\
\f{\pa E_{\omega}}{\pa q^1}=G_2(q)=\sum ^{b}_{\alpha=k+1}\mathcal{N}_2^{\,  \alpha}\widetilde{C}_{\alpha-k}\left(q\right)
\end{array}
\right. \ee and taking into account (\ref{j=1}), we obtain \ba
\hspace{-1cm}E_{\omega}(q^1,...,q^m,\dot{q}^1,...,\dot{q}^m)&=&
\sum_{a=1}^{k} \left[\int
\widetilde{L}_a\left(\dot{q}\right)\mathcal{N}^{\, a}_{i}
\dot{q}^i\mathcal{N}^ {a}_{1}d\dot{q}^1+\int
\widetilde{L}_a\left(\dot{q}\right)
 \mathcal{N}^{\,  a}_{i}\dot{q}^i\mathcal{N}^
{a}_{2}d\dot{q}^2
-\right.\nonumber\\
\hspace{-1cm}&& \int\int
\widetilde{L}_a'\left(\dot{q}\right)\mathcal{N}^{\,
a}_{i}\dot{q}^i \mathcal{N}^ {a}_{1}\mathcal{N}^
{a}_{2}d\dot{q}^1d\dot{q}^2-\nonumber\\
&&
\left.\int\int \widetilde{L}_a\left(\dot{q}\right)
\mathcal{N}^
{a}_{1}\mathcal{N}^
{a}_{2}d\dot{q}^1d\dot{q}^2\right]
+f_2(\dot{q}^3,...,\dot{q}^m)+\nonumber\\
\hspace{-1cm}&& \sum_{\alpha=k+1}^{b}\left[\int
\widetilde{C}_{\alpha-k}\left(q\right)\mathcal{N}_1^{\,
\alpha}dq^1+ \int
\widetilde{C}_{\alpha-k}\left(q\right)\mathcal{N}_2^{\,
\alpha}dq^2
-\right.\nonumber\\
&&\left.\int\int \widetilde{C}'_{\alpha-k}\left(q\right)\mathcal{N}_1^{\,  \alpha}
\mathcal{N}_2^{\,  \alpha}dq^1dq^2\right]
+
g_2(q^3,...,q^m)
\label{j=2}
\ea
which can be written in the form
\ba
E_{\omega}(q,\dot{q})&=&
\sum_{a=1}^{k}\sum_{l=1}^2\sum_{i_1<...<i_l=1}^{2}(-1)^{l+1}\underbrace{\int_{...} \int}_l
\left[
\widetilde{L}_a^{(l-1)}(\dot{q})\mathcal{N}^{ a}_{i}
\dot{q}^i+(l-1)\widetilde{L}_a^{(l-2)}(\dot{q})\right]\mathcal{N}^{ a}_{i_1}...\mathcal{N}^{ a}_{i_l}
d\dot{q}^{i_1}...d\dot{q}^{i_l}+\nonumber\\
&&\sum_{\alpha=k+1}^{b}\sum_{l=1}^2\sum_{i_1<...<i_l=1}^{2}(-1)^{l+1}\underbrace{\int_{...} \int}_l
\widetilde{C}_{\alpha-k}^{(l-1)}(q)\mathcal{N}^{\alpha}_{i_1}...\mathcal{N}^{\alpha}_{i_l}
dq^{i_1}...dq^{i_l}+\nonumber\\
&&f_2(\dot{q}^3,...,\dot{q}^m)+g_2(q^3,...,q^m)\nonumber\\
&&
\ea
where  $\widetilde{C}^{(l)}_{\alpha}:=\f {d^{l}\widetilde{C}_{\alpha}(\eta)}{d\eta^l}$, $\widetilde{L}^{(l)}_{a}:=\f 
{d^l\widetilde{L}_a(\eta)}{d\eta^l}$.

\noindent Repeating this procedure for j = 3, ... m, finally in
the $m-th$ and last step, we obtain \ba E_{\omega}(q,\dot{q})&=&
\sum_{a=1}^{k}\sum_{l=1}^m\sum_{i_1<...<i_l=1}^{m}(-1)^{l+1}\underbrace{\int_{...}
\int}_l \left[ \widetilde{L}_a^{(l-1)}(\dot{q})\mathcal{N}^{
a}_{i}
\dot{q}^i+(l-1)\widetilde{L}_a^{(l-2)}(\dot{q})\right]\mathcal{N}^{
a}_{i_1}...\mathcal{N}^{ a}_{i_l}
d\dot{q}^{i_1}...d\dot{q}^{i_l}+\nonumber\\
&&\sum_{\alpha=k+1}^{b}\sum_{l=1}^m\sum_{i_1<...<i_l=1}^{m}(-1)^{l+1}\underbrace{\int_{...}
\int}_l
\widetilde{C}_{\alpha-k}^{(l-1)}(q)\mathcal{N}^{\alpha}_{i_1}...\mathcal{N}^{\alpha}_{i_l}
dq^{i_1}...dq^{i_l}.\, \, \quad
\blacksquare
\label{energienonlinear}\ea

 Let us now discuss the question, what to do when the Birkhoffian given by (\ref{bir}) is not regular in the sense of definition 
(\ref{regular}).
\\
\textit{ If there exists at least one
loop in an LC circuit that contains only capacitors, then the Birkhoffian associated to the network is  \textbf{\textit{never} 
regular}}.
\\

\noindent Indeed, for the $l$-loop which contains only capacitors,
on the column $l$ of the matrix $A$  we have $A^a_{l}=0$ for any $a=1,...,k$. Without loss of  
generality, we will assume that $l=1$, that is  \be
 A^a_{1}=0, \quad \textrm{for any } a=1,...,k \label{a}
 \ee
For the Birkhoffian (\ref{bir}), the determinant in
(\ref{regular}) becomes \be
\textrm{det}\left[\f{\pa Q_j}{\pa \ddot{q}^i}(q,\, \dot{q},\, 
\ddot{q})\right]_{i,j=1,...,m}=\textrm{det}\left[\sum^{k}_{a=1}\mathcal{N}_j^{\,
a}\mathcal{N}^ {a}_{\,
i}\widetilde{L}_a\left(\dot{q}\right)\right]_{i,j=1,...,m} \ee
From (\ref{equivalence}), we get
$\mathcal{N}_j^a=\sum_{i_1=1}^{m}\mathcal{C}^{i_1}_jA^a_{i_1}$ for
any $a=1,...,k$, $j=1,...,m$.\\  Then,  taking into account
(\ref{a}), we have, for example, in the case $m=2$ \be
\sum^{k}_{a=1}\mathcal{N}_j^{\, a}\mathcal{N}^ {a}_{\,
i}\widetilde{L}_a\left(\dot{q}\right)=\mathcal{C}^2_j\mathcal{C}^2_i\left[
(A^1_{2})^2\widetilde{L}_1\left(\dot{q}\right)+(A^2_{2})^2\widetilde{L}_2\left(\dot{q}\right)+...+(A^k_{2})^2\widetilde{L}_k\left
(\dot{q}\right)\right] \ee Then, \be
\textrm{det}\left[\sum^{k}_{a=1}\mathcal{N}_j^{\,  a}\mathcal{N}^
{a}_{\,  i}\widetilde{L}_a\left(\dot{q}\right)\right]_{j,i=1,2}=
\left[\sum_{a=1}^k(A^a_{2})^2\widetilde{L}_a\left(\dot{q}\right)\right]^2
\left|
\begin{array}{lll}
\mathcal{C}^2_1\mathcal{C}^2_1& \mathcal{C}^2_1\mathcal{C}^2_2\\
\\
\mathcal{C}^2_1\mathcal{C}^2_2 &\mathcal{C}^2_2\mathcal{C}^2_2
\end{array}
\right|=0 \ee since the second factor obviously vanishes.  In the
case  $m=3$, we obtain \ba \sum^{k}_{a=1}\mathcal{N}_j^{\,
a}\mathcal{N}^ {a}_{\,
i}\widetilde{L}_a\left(\dot{q}\right)&=&\mathcal{C}^2_j\mathcal{C}^2_i\left[
\sum_{a=1}^k(A^a_{2})^2\widetilde{L}_a\left(\dot{q}\right)\right]
+\left(\mathcal{C}^2_j\mathcal{C}^3_i+\mathcal{C}^2_i\mathcal{C}^3_j\right)
\left[
\sum_{a=1}^kA^a_{2}A^a_{3}\widetilde{L}_a\left(\dot{q}\right)\right]+\nonumber\\
&&\mathcal{C}^3_j\mathcal{C}^3_i\left[
\sum_{a=1}^k(A^a_{3})^2\widetilde{L}_a\left(\dot{q}\right)\right]\label{elemente}
\ea Using basic calculus, the determinant of the matrix with
elements (\ref{elemente}) can be rearranged as a linear
combination of determinants having the columns of the form
$\left(\begin{array}{lll}
\mathcal{C}^{i_1}_1\mathcal{C}^{j_1}_1\\
\mathcal{C}^{i_1}_1\mathcal{C}^{j_1}_2\\
\mathcal{C}^{i_1}_1\mathcal{C}^{j_1}_3
\end{array}\right)$,
$\left(\begin{array}{lll}
\mathcal{C}^{i_1}_2\mathcal{C}^{j_1}_1\\
\mathcal{C}^{i_1}_2\mathcal{C}^{j_1}_2\\
\mathcal{C}^{i_1}_2\mathcal{C}^{j_1}_3
\end{array}\right)$,
$\left(\begin{array}{lll}
\mathcal{C}^{i_1}_3\mathcal{C}^{j_1}_1\\
\mathcal{C}^{i_1}_3\mathcal{C}^{j_1}_2\\
\mathcal{C}^{i_1}_3\mathcal{C}^{j_1}_3
\end{array}\right)$, respectively,  with $i_1, j_1=2$ or $3$ in each case. Hence, each of those determinants contain at least two 
linearly dependent columns, that is, they vanish, and this shows
that the determinant  is zero  in the case $m=3$ as well.
Similarly, for an arbitrary $m$, the determinant of the matrix
with the elements \ba \sum^{k}_{a=1}\mathcal{N}_j^{\,
a}\mathcal{N}^ {a}_{\,
i}\widetilde{L}_a\left(\dot{q}\right)&=&\sum_{i_1=2}^m\mathcal{C}^{i_1}_j\mathcal{C}^{i_1}_i\left[
\sum_{a=1}^k(A^a_{i_1})^2\widetilde{L}_a\left(\dot{q}\right)\right]
+\nonumber\\
&&\sum_{2<i_1<j_1}^m\left(\mathcal{C}^{i_1}_j\mathcal{C}^{j_1}_i+\mathcal{C}^{i_1}_i\mathcal{C}^{j_1}_j\right)
\left[
\sum_{a=1}^kA^a_{i_1}A^a_{j_1}\widetilde{L}_a\left(\dot{q}\right)\right]
\label{m}\ea  is zero.  $\quad \blacksquare$
\\

\noindent \textit{If
 there  exists in the network  $m'<m$  loops which contain only capacitors,
 all the other loops containing at least an inductor, we can
 \textbf{regularize} the Birkhoffian (\ref{bir}) via \textbf{reduction of the
 configuration space}. The reduced configuration space $\bar{M}_{c}$
 of dimension $m-m'$, is a linear subspace of $M_c$ or a manifold,
depending on  whether the capacitors  are linear or nonlinear. We
claim that  the Birkhoffian $\bar{\omega}_c$ of the reduced
configuration space $\bar{M}_c$ is  still a \textbf{conservative
Birkhoffian}. Under certain conditions on the functions $L_a$,
$a=1,...,k,$ which characterize the inductors, the reduced
Birkhoffian $\bar{\omega}_c$ will be a \textbf{regular
Birkhoffian}}.
\\

\noindent Without loss of generality, we can assume that there is
one loop in the network that contains only capacitors and in the
coordinate system we have chosen \be  \mathcal{N}_1^{a}=0, \,
a=1..., k\label{ngamma}\ee  Thus, the Birkhoffian components
(\ref{bir}), with (\ref{bir2}),  (\ref{bir1}), are given by,
$j=2,...,m$, \vspace{-0.5cm}\ba
 \hspace{-0.75cm}
 Q_1(q,\dot{q},\ddot{q})&=&\sum ^{b}_{\alpha=k+1}\mathcal{N}_1^{\alpha}\widetilde{C}_{\alpha-k}(q)\nonumber\\
\hspace{-0.75cm}Q_j(q,\dot{q},\ddot{q})&=&\sum^{m}_{i=2}
\sum^{k}_{a=1}\mathcal{N}^a_{j}\mathcal{N}^
{a}_{i}\widetilde{L}_{a}\left(\dot{q}\right)\ddot{q}^i+\sum
^{b}_{\alpha=k+1}\mathcal{N}_j^{\alpha}\widetilde{C}_{\alpha-k}(q)
\ea We note that, according to (\ref{ngamma}),
 $\dot{q}^1$ does  not
appear in any function $\widetilde{L}_{a}(\dot{q})$
 and  the terms
$\widetilde{L}_{a}(\dot{q})\ddot{q}^1$ do not appear in any of the
Birkhoffian components $Q_2(q,\dot{q}, \ddot{q}),..., Q_m(q,
\dot{q}, \ddot{q})$.\\
If the capacitors in this loop are linear devices, $Q_1$ is a
linear combination of $q$'s and
we can use this relation to reduce the  configuration 
space  $M_c$, to an affine-linear subspace   $\bar{M}_c$ of
dimension $m-1$. If the capacitors in this loop are nonlinear
devices, $Q_1$ depends nonlinearly on the $q$'s. We define the
$(m-1)$-dimensional manifold  $\bar{M}_{c}\subset M_c$ by \be
\bar{M}_c=\{q\in M_c\,\, | \, \, \sum
^{b}_{\alpha=k+1}\mathcal{N}_1^{\alpha}\widetilde{C}_{\alpha-k}(q)=0\}
\ee
 By the implicit function theorem, we obtain a
local coordinate system on the reduced configuration space
$\bar{M}_{c}$. Taking $\bar{q}^1:=q^2$,..., $\bar{q}^{m-1}:=q^m$,
the Birkhoffian has the form
$\bar{\omega}_c=\sum^{m-1}_{j=1}\bar{Q}_jd\bar{q}^j$,  \be
\bar{Q}_j(\bar{q}, \dot{\bar{q}},
\ddot{\bar{q}})=\bar{F}_j(\dot{\bar{q}})\ddot{\bar{q}}+\bar{G}_j(\bar{q}),
\quad \textrm{where} \label{birred} \ee
 \be
\bar{F}_j(\dot{\bar{q}})\ddot{\bar{q}}:= \sum^{m-1}_{i=1}
\sum^{k}_{a=1}\mathcal{N}^a_{(j+1)}\mathcal{N}^ {a}_{(i+1)}
L_{a}\left(\sum^{m-1}_{l=1} \mathcal{N}^a_{(l+1)}
\dot{\bar{q}}^l\right) \ddot{\bar{q}}^i\ee
 \be
\bar{G}_j(\bar{q}):= \sum
^{b}_{\alpha=k+1}\mathcal{N}_{(j+1)}^{\alpha} C_{\alpha-k}
\left(\mathcal{N}^\alpha_1f(\bar{q}^1,...,\bar{q}^{m-1})+\sum^{m-1}_{l=1}
\mathcal{N}^\alpha_{(l+1)} \bar{q}^l+const\right)
 \label{bir6'} \ee
 $f:U\subset\mathbf{R}^{m-1}\longrightarrow \mathbf{R}$ being the unique function such that
$f(\bar{q}_0)=q^1_0$, $q^1_0\in \mathbf{R}$, and \be \sum
^{b}_{\alpha=k+1}\mathcal{N}_{1}^{\alpha} C_{\alpha-k}
\left(\mathcal{N}^\alpha_1f(\bar{q}^1,...,\bar{q}^{m-1})+\sum^{m-1}_{l=1}
\mathcal{N}^\alpha_{(l+1)}
\bar{q}^l+const\right)=0\label{implicit} \ee
 for all $\bar{q}=(\bar{q}^1,..., \bar{q}^{m-1})\in U$, with
$U$ a neighborhood of
$\bar{q}_0=(\bar{q}_0^1,...,\bar{q}_0^{m-1})$.\\
 We will now prove that the Birkhoffian (\ref{birred}) is
conservative.
 In order to do so, we will show that there exists  a
function $\bar{E}_{\omega}(\bar{q},\dot{\bar{q}})$ satisfying \be
\sum^{m-1}_{j=1}\bar{Q}_j(\bar{q}\, \, \dot{\bar{q}},\,
\ddot{\bar{q}})\dot{\bar{q}}^j=\sum^{m-1}_{j=1}\left[\f{\pa
\bar{E}_{\omega}}{\pa \bar{q}^j}\dot{\bar{q}}^j+\f{\pa
\bar{E}_{\omega}}{\pa \dot{\bar{q}}^j}\ddot{\bar{q}}^j\right]
\label{dissip'} \ee Because of the special form of the terms on
the left side of (\ref{dissip'}), we may assume that
$\bar{E}_{\omega}(\bar{q},\dot{\bar{q}})$ is a sum of a function
 depending only on $\bar{q}$, and a function  depending only on
$\dot{\bar{q}}$. From the theory of total differentials, a
necessary condition for the existence of such functions is the
fulfillment of the following relations \be \left\{
\begin{array}{ll}
\f {\partial \bar{\mathcal{F}}_j(\dot{\bar{q}})}{\partial
\dot{\bar{q}}^l }-\f {\partial \bar{\mathcal{F}}_l
(\dot{\bar{q}})}{\partial
\dot{\bar{q}}^j}=0\\
\\
\f {\partial \bar{G}_j(\bar{q})}{\partial \bar{q}^l}-\f {\partial
\bar{G}_l(\bar{q})}{\partial \bar{q}^j}=0
\end{array}
\right.\label{bir5'} \ee for any $j,l=1,...,m-1$, where \be
\bar{\mathcal{F}}_j(\dot{\bar{q}}):=\sum^{m-1}_{i=1}\sum^{k}_{a=r+1}\mathcal{N}^{\,
a}_{(j+1)}\mathcal{N}^ {a}_{(i+1)} L_{a} \left(\sum^{m-1}_{l=1}
\mathcal{N}^a_{(l+1)} \dot{\bar{q}}^l\right) \dot{\bar{q}}^i
 \ee    We  check in the same way as for the functions $\mathcal{F}_j(\dot{q})$ in
 (\ref{bir6}), that the first relation in (\ref{bir5'})
is fulfilled. From (\ref{bir6'}), the second relation in
(\ref{bir5'}) reads as \ba &&\sum
^{b}_{\alpha=k+1}\left\{\mathcal{N}_{(j+1)}^{\alpha}
\widetilde{C}'_{\alpha-k}(\bar{q})\left[\mathcal{N}^\alpha_1\f{\pa
f(\bar{q})}{\pa \bar{q}^l}
+\mathcal{N}_{(l+1)}^{\alpha}\right]-\right.\nonumber\\
&& \left. \quad \quad \quad \quad\mathcal{N}_{(l+1)}^{\alpha}
\widetilde{C}'_{\alpha-k}(\bar{q})\left[\mathcal{N}^\alpha_1\f{\pa
f(\bar{q})}{\pa \bar{q}^j}
+\mathcal{N}_{(j+1)}^{\alpha}\right]\right\}=0 \label{c}\ea where
$\widetilde{C}'_{\alpha-k}:=\f
{d\widetilde{C}_{\alpha-k}(\eta)}{d\eta}$. The relation (\ref{c})
reduces to \be \hspace{-0.7cm}\sum ^{b}_{\alpha=k+1}
\mathcal{N}_{(j+1)}^{\alpha}
\widetilde{C}'_{\alpha-k}(\bar{q})\mathcal{N}^\alpha_1\f{\pa
f(\bar{q})}{\pa \bar{q}^l} - \mathcal{N}_{(l+1)}^{\alpha}
\widetilde{C}'_{\alpha-k}(\bar{q})\mathcal{N}^\alpha_1\f{\pa
f(\bar{q})}{\pa \bar{q}^j}  =0 \label{nec2} \ee Taking into
account  (\ref{implicit}), the above relation is fulfilled, for
any $j,l=1,...,m-1$. Indeed, taking the derivatives with respect
to $\bar{q}^j$ and also to $\bar{q}^l$, in the equation
(\ref{implicit}), we obtain, respectively, \ba \sum
^{b}_{\alpha=k+1}\mathcal{N}^\alpha_1\widetilde{C}'_{\alpha-k}(\bar{q})
\left[\mathcal{N}^\alpha_1\f{\pa f(\bar{q})}{\pa \bar{q}^j}
+\mathcal{N}_{(j+1)}^{\alpha}\right]&=&0\nonumber\\
\sum
^{b}_{\alpha=k+1}\mathcal{N}^\alpha_1\widetilde{C}'_{\alpha-k}(\bar{q})
\left[\mathcal{N}^\alpha_1\f{\pa f(\bar{q})}{\pa \bar{q}^l}
+\mathcal{N}_{(l+1)}^{\alpha}\right]&=&0\label{nec3} \ea
 Now we
multiply  in (\ref{nec3}) the first equation  with $\f{\pa
f(\bar{q})}{\pa \bar{q}^l}$, the second  equation with $-\f{\pa
f(\bar{q})}{\pa \bar{q}^j}$ and we add the resulting  equations to
obtain the equation (\ref{nec2}).\\
Thus, we proved  the existence of  a function
$\bar{E}_{\omega}(\bar{q},\dot{\bar{q}})$ such that
(\ref{dissip'}) is fulfilled.
\\

\noindent For any other  loop which  contains only capacitors, we just repeat this procedure. Thus, we finally arrive at a configuration space 
$\bar{M}_c$ of dimension $m-m'$, where $m'$ denotes the total
number of  loops of that type.  $\quad \blacksquare$
\\

In case the network has loops  which contain only inductors  the
Birkhoffian can be a regular one but we can further reduce the
configuration space. Inductor loops can be considered as some
conserved quantities of the network.
\\

\noindent \textit{If
 there  exists in the network  $m''<m$  loops which contain only linear inductors,
 all the other loops containing at least a capacitor, we can
 further reduce the configuration space.  The reduced configuration space $\hat{M}_{c}$
 of dimension $m-m''$, is a linear  subspace of $M_c$. We
claim that  the Birkhoffian $\hat{\omega}_c$ of the reduced
configuration space $\hat{M}_c$ is  a \textbf{conservative
Birkhoffian}.  Under certain conditions on the functions $L_a$,
$a=1,...,k,$ which characterize the inductors, the reduced
Birkhoffian $\hat{\omega}_c$ will be a \textbf{regular
Birkhoffian}. }

 \noindent Without loss of generality, we can
assume that there is one loop in the network that contains only
inductors and in the coordinate system we have chosen \be
\mathcal{N}_1^{\alpha}=0, \, \alpha=1..., p\label{nalpha}\ee Thus,
the Birkhoffian components (\ref{bir}), with (\ref{bir2}),
(\ref{bir1}), are given by, $j=2,...,m$, \ba
 \hspace{-0.75cm}
 Q_1(q,\dot{q},\ddot{q})&=&\sum^{m}_{i=1}
\sum^{k}_{a=1}\mathcal{N}^a_{1}\mathcal{N}^
{a}_{i}\widetilde{L}_{a}\left(\dot{q}\right)\ddot{q}^i\nonumber\\
\hspace{-0.75cm}Q_j(q,\dot{q},\ddot{q})&=&\sum^{m}_{i=1}
\sum^{k}_{a=1}\mathcal{N}^a_{j}\mathcal{N}^
{a}_{i}\widetilde{L}_{a}\left(\dot{q}\right)\ddot{q}^i+\sum
^{b}_{\alpha=k+1}\mathcal{N}_j^{\alpha}\widetilde{C}_{\alpha-k}(q)
\ea We note that, according to (\ref{nalpha}),
 $q^1$ does  not
appear in any function $\widetilde{C}_{\alpha-k}(q)$.\\
If the inductors in this loop are linear devices, $Q_1$ is a
linear combination of $\ddot{q}$'s. We can integrate this relation
to obtain an affine-linear relation between $q$'s (see the first
example in section 5).
We can use this relation to reduce the  configuration 
space  $M_c$, to an affine-linear subspace   $\hat{M}_c$ of
dimension $m-1$. Taking $\hat{q}^1:=q^2$,...,
$\hat{q}^{m-1}:=q^m$, one can write the Birkhoffian components of
$\hat{\omega}_{c}$ and one can prove, using the same ideas as in
the previous reduction case, the existence of the function $
\hat{E}_{\omega}$ such that this Birkhoffian is conservative.
\\

\noindent For any other  loop which  contains only linear inductors,
we just repeat this procedure. Thus, we finally arrive at a configuration space 
$\hat{M}_c$ of dimension $m-m''$, where $m''$ denotes the total
number of  loops of that type. $\quad \blacksquare$
\\

If the devices in the $m"<m$ inductor loops are nonlinear devices,
then, $Q_1(q,\dot{q},\ddot{q})$,..., $Q_{m"}(q,\dot{q},\ddot{q})$,
are nonlinear functions depending on $\dot{q}$'s and $\ddot{q}$'s.
Using these relations we can define a smooth constant rank affine
sub-bundle
$\mathfrak{S}_c$ of the affine bundle $\pi_J:J^2(M_c)\longrightarrow TM_c$, on which we define the constrained 
Birkhoffian system ($M_c, \omega_c, \mathfrak{S}_c$).
 The submanifold
$\mathfrak{S}_c$ has codimension $m''$.

\section{LC electric circuits with independent current/voltage
sources}

\noindent Let us now consider  an electric  circuit containing {\scriptsize S}$_I$ independent current sources and {\scriptsize S}$_V$ 
independent voltage sources, in addition to  $k$ inductors and $p$ capacitors.
Then $b$$=k+p$+{\scriptsize S}$_I$+{\scriptsize 
S}$_V$=$m+n$, where $b$, $m$, $n$ have the same meaning as in
section 3. We suppose that  $m-{\textrm{\scriptsize S}}_I>0$,
$n-{\textrm{\scriptsize S}}_V>0$. The branches of the oriented
connected graph associated to this  circuit are labelled as
follows: $\textsc{l}_a$, $a=1,...,k$, the inductor branches,
$\textsc{c}_\alpha$, $\alpha=1,...,p$, the capacitor branches,
$S_{I_i}$, $i=1,...,${\scriptsize S}$_I$, the current source
branches, and $S_{V_j}$,
$j=1,..., ${\scriptsize 
S}$_V$, the voltage source branches.

\noindent Let the basic equations governing the circuit be now
written in the form \be \left\{
\begin{array}{ll}
\mathcal{B}_1^T
\left(
\begin{array}{c}
\textsc{i}_a \\
\f{d\textsc{q}_{\alpha}}{dt}
\end{array}
\right)+\mathcal{B}_2^T
\left(
\begin{array}{c}
\textsc{i}_{s_I}(t)
\end{array}
\right)=0\\
\\
\mathcal{A}_1^T
\left(
\begin{array}{c}
L_{a}(\textsc{i}_a)\,
\f{d\textsc{i}_a}{dt} \\
C_\alpha(\textsc{q}_{\alpha})
\end{array}
\right)
+\mathcal{A}_2^T
\left(
\begin{array}{c}
v_{s_V}(t)
\end{array}
\right)=0
\end{array}
\right.\label{surse7*}
\ee
where $ \mathcal{A}_1^T\in \mathfrak{M}_{(m-{\textrm{\scriptsize S}_I})(k+p)}(\mathbf{R})$, $\mathcal{A}_2^T\in  
\mathfrak{M}_{(m-
{\textrm{\scriptsize S}_I}){\textrm{\scriptsize S}}_V}(\mathbf{R})$,
 $\mathcal{B}_1^T\in  \mathfrak{M}_{(n-{\textrm{\scriptsize S}_V})(k+p)}(\mathbf{R})$, $\mathcal{B}_2^T\in  
\mathfrak{M}_{(n-{\textrm{\scriptsize S}_V}){\textrm{\scriptsize
S}}_I}(\mathbf{R})$. We also assume that
rank$(\mathcal{A}_1^T)=m-{\textrm{\scriptsize S}_I}$,
$(\mathcal{B}_1^T)=n-{\textrm{\scriptsize S}_V}$. The functions
$\textsc{i}_{s_I}(t)$ and  $v_{s_V}(t)$  are given vector
functions of time. They describe  the independent current sources
and independent voltage sources,  respectively.
 The other quantities in (\ref{surse7*}) are  defined as in section 3.
\\

\noindent \textit{In the following we give a Birkhoffian
formulation for the network described by the system of equations
(\ref{surse7*}), using the same procedure as in section 3. That
is,
 using the first set of equations (\ref{surse7*}),
 we are going to define a  family of
  ($m-\textrm{\scriptsize S}_I$)-dimensional affine-linear configuration spaces
  $M_c\subset\mathbf{R}^b$ parameterized by a constant vector
  $c$ in $\mathbf{R}^{n-{\textrm{\scriptsize S}_V}}$.
A Birkhoffian $\omega_{t_c}$ on the configuration space $M_c$
arises from a linear combination of the second set of equations
(\ref{surse7*}). Thus, ($M_c,\omega_{t_c}$) will be a family of
Birkhoff systems that describe the  LC circuit with independent
current/voltage sources considered.}
\\

\noindent Let $\mathfrak{H}: \mathbf{R}^{k+p}
\longrightarrow\mathbf{R}^{n -{\textrm{\scriptsize S}_V}}$ be the
linear map  that, with respect to a coordinate system
($x^1,...,x^{k+p}$) on $\mathbf{R}^{k+p}$, is given by \be
\mathfrak{H}(x^1,...,x^{k+p})=\mathcal{B}_1^T \left(
\begin{array}{c}
x^1 \\
\vdots\\
x^{k+p}
\end{array}
\right)+\mathcal{B}_2^T
\left(
\begin{array}{c}
\textsc{i}_{s_I}(t)
\end{array}
\right)\label{surse9} \ee We define \be
M_c:=\mathfrak{H}^{-1}(c)\label{surse8'} \ee $c$ being a constant
vector in $\mathbf{R}^{n-{\textrm{\scriptsize S}_V}}$. $M_c$ is a
time-dependent affine linear subspace in $\mathbf{R}^{k+p}$. From
rank$(\mathcal{B}_1^T)=n-{\textrm{\scriptsize S}_V} $, its
dimension is $ k+p+\textrm{\scriptsize
S}_V-n=m-\textrm{\scriptsize S}_I$.

\noindent Let us figure out  the relation between $\textsc{i}_a$,
$\f{d\textsc{q}_{\alpha}}{dt}$, and  coordinates  on $M_c$. As in
the case without sources, taking into account (\ref{var}) and the
fact that the matrix $\mathcal{B}_1^T$ is a constant matrix, we
integrate  the first set of equations (\ref{surse7*}) to arrive at
\be \mathcal{B}_1^T \left(
\begin{array}{c}
\textsc{q}_{(a)} \\
\textsc{q}_{\alpha}
\end{array}
\right)+\mathcal{I}(t)=c\label{surse7''} \ee with $c$ a constant
vector in $\mathbf{R}^{n-{\textrm{\scriptsize S}_V}}$ and
$\mathcal{I}(t)$ a primitive of $\mathcal{B}_2^T \left(
\begin{array}{c}
\textsc{i}_{s_I}(t)
\end{array}
\right)$. \\
Likewise  consider coordinates  in $\mathbf{R}^{k+p}$ \be
x^1:=\textsc{q}_{(1)},\,
..,\,x^k:=\textsc{q}_{(k)},\,x^{k+1}:=\textsc{q}_{1},\,..,\,
x^{k+p}:=\textsc{q}_{p} \label{cssurse} \ee We can define
coordinates on $M_c$ by solving the equations (\ref{surse7''}) in
terms of an appropriate set of $(m-\textrm{\scriptsize S}_I)$ of
the $\textsc{q}$-variables, say
$q=(q^1,...,q^{m-\textrm{\scriptsize S}_I})$. In other words, we
express any of the $x$-variables as a function of
$q=(q^1,...,q^{m-\textrm{\scriptsize S}_I})$, namely, \ba &&x^a=
\sum_{j=1}^{m-\textrm{\scriptsize S}_I}\mathfrak{N}^{a}_jq^j+f^a(t)+const,\, a=1,...,k,\nonumber\\
&&x^{\alpha}= \sum_{j=1}^{m-\textrm{\scriptsize
S}_I}\mathfrak{N}^{\alpha}_jq^j+f^{\alpha}(t)+const, \,
\alpha=k+1,...,k+p \label{surse10'} \ea with certain constants
$\mathfrak{N}^{a}_{ j}$, $\mathfrak{N}^{\alpha}_{ j}$ and certain
functions of $t$, $f^a(t)$, $f^{\alpha}(t)$. \\
The constant matrix
$\mathfrak{N}=\left(\begin{array}{c}
\mathfrak{N}^{a}_{ j}\\
\mathfrak{N}^{\alpha}_{j}
\end{array}\right)_{{a=1,...,k,\,  \alpha=k+1,...,k+p \atop j=1,...,m-\textrm{\scriptsize S}_I}}$ has rank $m-\textrm{\scriptsize 
S}_I$,   and  there exists a nonsingular matrix $\mathcal{C}\in \mathfrak{M}_{(m-\textrm{\scriptsize S}_I)(m-\textrm{\scriptsize 
S}_I)}(\mathbf{R})$   such that \be
\mathcal{C}\mathcal{A}_1^T=\mathfrak{N}^T\label{sequivalence} \ee

\noindent \textit{We define the Birkhoffian $\omega_{t_c}$ of
$M_c$ such that the differential system (\ref{difsistem}) is a
linear combination of the second set of equations in
(\ref{surse7*}), which is obtained multiplying the second set of
equations in
 (\ref{surse7*}) by the matrix
$\mathcal{C}$. Taking into account (\ref{sequivalence}), in terms
of q-coordinates as chosen before, the expressions of the
components $Q_j(t,q, \, \dot{q},\,\ddot{q})$,
$j=1,...,m-\textrm{\scriptsize S}_I $ are} \ba Q_j(t,q, \,
\dot{q},\,\ddot{q})&=&F_j(t,\dot{q})\ddot{q}+G_j(t,q)+\mathcal{V}_{j}(t)\label{sursebir}
\ea \textit{where }\ba F_j(t,\dot{q})\ddot{q}&=&
\sum^{k}_{a=1}\mathfrak{N}^a_{j}\widetilde{L}_a
\left(t,\dot{q}\right) \left(\sum^{m-\textrm{\scriptsize
S}_I}_{i=1}\mathfrak{N}^ {a}_{i}\ddot{q}^i+\f{d^2f^a(t)}{dt^2}
\right)=\nonumber\\
&&
\sum^{m-\textrm{\scriptsize S}_I}_{i=1}
\left(
\sum^{k}_{a=1}\mathfrak{N}^a_{j}\mathfrak{N}^
{a}_{i}\widetilde{L}_a
\left(t,\dot{q}\right)
\right)\ddot{q}^i+\sum^{k}_{a=1}\mathfrak{N}^a_{j}\widetilde{L}_a\left(t,\dot{q}\right)\f{d^2f^a(t)}{dt^2}
\label{surse13}
\ea
\ba
G_j(t,q)
=\sum ^{k+p}_{\alpha=k+1}\mathfrak{N}_j^{\alpha}C_{\alpha-k}\left(\sum_{j=1}^{m-\textrm{\scriptsize 
S}_I}\mathfrak{N}^{\alpha}_jq^j+f^{\alpha}(t)+const\right)
=\sum ^{k+p}_{\alpha=k+1}\mathfrak{N}_j^{\alpha}\widetilde{C}_{\alpha-k}\left(t,q\right)
\label{surse11}
\ea
\ba
\mathcal{V}_{j}(t)=\sum^{\textrm{\scriptsize B}}_
              {s_V=k+p+\textrm{\scriptsize S}_I+1}
(\mathcal{C}^T\mathcal{A}_2^T)_{js_V}v_{s_V-k-p-S_I}(t)
\label{surse14} \ea
\\

\noindent \textit{If there exist in the network
$m'<m-\textrm{\scriptsize S}_I$ loops which contain only
capacitors or capacitors and independent voltage sources, then the
Birkhoffian
associated to the network  is {\bf \textit{never} regular}}.
\\

\noindent Indeed,  in this case  the  functions $Q_j$ corresponding to such   loops 
depend only on $q$'s and $t$. Using the same procedure as in
section 3, the reduced configuration space $\bar{M}_{c}$
 of dimension $(m-\textrm{\scriptsize S}_I)-m'$, is a linear subspace of $M_c$ or a manifold,
depending on  whether the devices in the loops  are linear or
nonlinear.
\\

Let us finally discuss the question whether the Birkhoffian (\ref{sursebir})  is  conservative or not.
\\

\noindent \textit{ For  a linear LC circuit with independent
current/voltage sources we claim that the Birkhoffian
(\ref{sursebir}) is conservative. A nonlinear LC circuit with
independent current/voltage sources is conservative
if and only if it does not contain cutsets of inductors and  independent 
current sources.}
\\

\noindent In order to show that the Birkhoffian (\ref{sursebir})
is conservative, we are looking for a smooth function
$E_{\omega_{t}}(t,q,\dot{q})$ such that the relation
(\ref{surseconserv}) is fulfilled. For the Birkhoffian
(\ref{sursebir}), this relation becomes \ba
\left(\sum^{m-\textrm{\scriptsize
S}_I}_{i=1}\sum^{k}_{a=1}\mathfrak{N}^a_{i}\mathfrak{N}^a_{j}\widetilde{L}_a
\left(t,\dot{q}\right)\dot{q}^i\right)\ddot{q}^j+
\left(\sum^{k}_{a=1}\mathfrak{N}^a_{j}\widetilde{L}_a
\left(t,\dot{q}\right)\f
{d^2f^a(t)}{dt^2}+G_j(t,q)+\mathcal{V}_{j}(t)\right)
\dot{q}^j=\nonumber \ea \ba \f{\pa E_{\omega_{t}}}{\pa
q^j}\dot{q}^j+ \f{\pa E_{\omega_{t}}}{\pa
\dot{q}^j}\ddot{q}^j\label{sfarsit} \ea

\noindent If the inductors and the capacitors in the network are
linear devices, taking into account  (\ref{15}), we easily find
the function \ba
E_{\omega_t}(t,q,\dot{q})&=&\f{1}{2}\sum^{k}_{a=1}\sum_{i,j=1}^{m-\textrm{\scriptsize
S}_I} {\textrm {\scriptsize L}}_a
\mathcal{N}^{a}_i\mathcal{N}^{a}_j\dot{q}^i \dot{q}^j+\f{1}{2}
\sum^{{\textrm{\scriptsize
B}}}_{\alpha=k+1}\sum^{m-\textrm{\scriptsize S}_I}_{i,j=1}
\f{\mathcal{N}^{\alpha}_i\mathcal{N}^{\alpha}_j}
{{\textrm {\scriptsize C}}_{\alpha-k}}q^iq^j+\nonumber\\
&&\sum_{j=1}^{m-\textrm{\scriptsize
S}_I}[\mathfrak{N}^a_{j}{\textrm {\scriptsize L}}_a \f
{d^2f^a(t)}{dt^2}+\mathcal{V}_{j}(t)+const_j]q^j
\label{energiesurse1} \ea which satisfies (\ref{sfarsit}).

\noindent If  the inductors and the capacitors in the network are
nonlinear devices, the existence of the function
$E_{\omega_t}(t,q,\dot{q})$ which satisfies (\ref{sfarsit}),
depends on the appearance of the  term
$\sum^k_{a=1}\mathfrak{N}^a_{j}\widetilde{L}_a
\left(t,\dot{q}\right)\f {d^2f^a(t)}{dt^2}$  in (\ref{sfarsit}).
For the networks which do not contain cutsets of inductors and  independent 
current sources, this term  does not appear at all in
(\ref{sfarsit}). In this case the proof of the existence of the
function $ E_{\omega_{t}}$  is the same as  in the case without
sources. If the term
$\sum^k_{a=1}\mathfrak{N}^a_{j}\widetilde{L}_a
\left(t,\dot{q}\right)\f {d^2f^a(t)}{dt^2}$ is different from zero
in (\ref{sfarsit}), then, $\f{\pa^2 E_{\omega_{t}}}{\pa q^j\pa 
\dot{q}^j}\neq \f{\pa^2 E_{\omega_{t}}}{\pa \dot{q}^j\pa q^j}$.
Therefore, the Birkhoffian (\ref{sursebir}) is not conservative in
the sense of definition (\ref{surseconserv}).
  $
\blacksquare$

\section{Examples}

The first example that we present is the example from the paper (\cite{maschke}), in which we have interchanged the capacitor $C_3$ 
and the inductor $L_1$ to emphasis that networks which contain capacitor loops and inductor cutsets fit into the formalism presented 
in section 3. The directed connected graph associated to this
circuit is presented in figure 1, page ~\pageref{figura}.

We first suppose that all devices are linear, that is,  they are
described by the relations (\ref{15}). Then,  taking into account
the values of the matrices $A$, $B$ given by (\ref{ex91}), the
equations (\ref{7}) which govern the network have the form \be
\left\{
\begin{array}{llllllll}
\textsc{i}_4-\f{d\textsc{q}_1}{dt}+\f{d\textsc{q}_3}{dt}=0\\
\textsc{i}_2+\textsc{i}_3-\f{d\textsc{q}_2}{dt}-\f{d\textsc{q}_3}{dt}=0\\
\textsc{i}_1-\textsc{i}_3-\textsc{i}_4=0\\$$
\\
\f {\textsc{q}_1}{\textrm{\scriptsize C}_1}-\f {\textsc{q}_2}{\textrm{\scriptsize C}_2}+\f {\textsc{q}_3}{\textrm{\scriptsize 
C}_3}=0\\
\textrm{\scriptsize L}_2\f{d\textsc{i}_2}{dt}+\f {\textsc{q}_2}{\textrm{\scriptsize C}_2}=0\\
\textrm{\scriptsize L}_1\f{d\textsc{i}_1}{dt}-\textrm{\scriptsize L}_2\f{d\textsc{i}_2}{dt}+\textrm{\scriptsize L}_3\f{d\textsc{i}_3}{dt}=0\\
-\textrm{\scriptsize L}_1\f{d\textsc{i}_1}{dt}-\textrm{\scriptsize
L}_4\f{d\textsc{i}_4}{dt}-\f {\textsc{q}_1}{\textrm{\scriptsize
C}_1}=0
\end{array}
\right.\label{ex92}
\ee
where $\textrm{\scriptsize C}_\alpha\neq 0$, $\alpha=1,2,3$ and $\textrm{\scriptsize L}_a\neq 0$, $a=1,2,3,4$, are distinct 
constants. The relations (\ref{var}), (\ref{cs}) read as follows
for this example \be \textsc{i}_a:=\f{d\textsc{q}_{(a)}}{dt},
\quad a=1,2,3,4\label{var9} \ee \be x^1:=\textsc{q}_{(1)},\,
...,\,x^4:=\textsc{q}_{(4)},\,x^{5}:=\textsc{q}_{1},...,\,
x^{7}:=\textsc{q}_{3}\label{cs9} \ee Using the first set of
equations (\ref{ex92}), we  define the 4-dimensional affine-linear
configuration  space $M_c$. We solve  the corresponding equations
(\ref{7''})  in terms of 4 variables. In view of  the notations
(\ref{var9}), (\ref{cs9}), we obtain, for example, \ba
x^{1}&=&x^3+x^{4}+const\nonumber\\
x^5&=&x^4+x^7+const\nonumber\\
x^{6}&=&x^{2}+x^{3}-x^{7}+const \label{ex94} \ea Thus, a
coordinate system on $M_c$ is given by \be q^1:=x^7, q^2:=x^2,
q^3:=x^3,q^4:=x^4\label{csex9} \ee The matrices of constants
$\mathcal{N}=\left(\begin{array}{c}
\mathcal{N}^{a}_{ j}\\
\mathcal{N}^{\alpha}_{j}
\end{array}\right)_{{a=1,2,3,4, \, \alpha=5,6,7 \atop j=1,2,3,4}}$  and  $\mathcal{C}$ in
 (\ref{10'}), (\ref{equivalence}) attain the form
\be
\mathcal{N}=\left(
\begin{array}{cccc}
 0& 0& 1& 1\\
 0& 1& 0& 0\\
 0& 0& 1& 0\\
 0& 0& 0& 1\\
 1& 0& 0& 1\\
 -1& 1& 1& 0\\
 1& 0& 0& 0
\end{array}
\right), \quad \mathcal{C}=\left(
\begin{array}{cccc}
 1& 0& 0& 0\\
 0& 1& 0& 0\\
 0& 1& 1& 0\\
 0& 0& 0& -1
\end{array}
\right)
\label{N9}
\ee
Note that, if we define the Birkhoffian $\omega_c$ of $M_c$  using  the second set of equations  (\ref{ex92}) and not a linear 
combination of them, that is,  the matrix $A^T$ instead of
$\mathcal{C}A^T=\mathcal{N}^T$, then, in terms of  the
$q$-coordinates introduced in (\ref{csex9}), we
 obtain
$\omega_c=\sum^{4}_{j=1}Q_j(q,\dot{q},\ddot{q})dq^j$, with \ba
&&Q_1(q,\dot{q},\ddot{q})=
\left(\frac{1}{\textrm{\scriptsize C}_1}+\f {1}{\textrm{\scriptsize C}_2}+\f {1}{\textrm{\scriptsize 
C}_3}\right)q^1-\f{q^2}{\textrm{\scriptsize C}_2}-\f{q^3}{\textrm{\scriptsize C}_2}+\f{q^4}{\textrm{\scriptsize 
C}_1}+const\nonumber\\
&&Q_2(q,\dot{q},\ddot{q})=\textrm{\scriptsize L}_2\ddot{q}^2-\f{q^1}{\textrm{\scriptsize C}_2}+\f{q^2}{\textrm{\scriptsize 
C}_2}+\f{q^3}{\textrm{\scriptsize C}_2}+const\nonumber\\
&&Q_3(q,\dot{q},\ddot{q})=-\textrm{\scriptsize L}_2\ddot{q}^2+(\textrm{\scriptsize L}_1+\textrm{\scriptsize L}_3)\ddot{q}^3+
\textrm{\scriptsize L}_1\ddot{q}^4\nonumber\\
&&Q_4(q,\dot{q},\ddot{q})=-\textrm{\scriptsize L}_1\ddot{q}^3-(\textrm{\scriptsize L}_1+\textrm{\scriptsize L}_4)
\ddot{q}^4-\f 
{q^1}{\textrm{\scriptsize C}_1}-\f{q^4}{\textrm{\scriptsize
C}_1}+const\label{alegere} \ea The  Birkhoffian (\ref{alegere}) is
\textbf{ \textit{not}}  conservative. Indeed, for the Birkhoffian
(\ref{alegere}) two of the necessary conditions for the existence
of the function $E_{\omega}:TM\to \mathbf{R}$ such that
(\ref{conserv}) is fulfilled, are \be \left\{
\begin{array}{lll}
\f{\pa E_{\omega}}{\pa \dot{q}^2}=\textrm{\scriptsize L}_2\dot{q}^2-\textrm{\scriptsize L}_2\dot{q}^3\\
\\
\f{\pa E_{\omega}}{\pa \dot{q}^3}=(\textrm{\scriptsize L}_1+\textrm{\scriptsize L}_3)\dot{q}^3+
\textrm{\scriptsize L}_1\dot{q}^4\end{array}
\right.
\ee
Because  $\textrm{\scriptsize L}_2\ne 0$, we see that $\f{\partial ^2 E_{\omega}}{\partial \dot{q}^3\dot{q}^2}\ne \f{\partial^2 
E_{\omega}}{\partial \dot{q}^2\dot{q}^3}$. Therefore, there does
not exist a function $E_{\omega}$ such that (\ref{conserv}) is
fulfilled.

However, proceeding as suggested in section 3,  the functions
$Q_j(q,\dot{q},\ddot{q})$, $j=1,2,3,4$ are given by (\ref{bir}),
(\ref{bir2}), and (\ref{bir1}), that is, \ba
&&Q_1(q,\dot{q},\ddot{q})=\left(\frac{1}{\textrm{\scriptsize C}_1}+\f {1}{\textrm{\scriptsize C}_2}+\f {1}{\textrm{\scriptsize 
C}_3}\right)q^1-\f{q^2}{\textrm{\scriptsize C}_2}-\f{q^3}{\textrm{\scriptsize C}_2}+\f{q^4}{\textrm{\scriptsize 
C}_1}+const\nonumber\\
&&Q_2(q,\dot{q},\ddot{q})=\textrm{\scriptsize L}_2\ddot{q}^2-\f{q^1}{\textrm{\scriptsize C}_2}+\f{q^2}{\textrm{\scriptsize 
C}_2}+\f{q^3}{\textrm{\scriptsize C}_2}+const\nonumber\\
&&Q_3(q,\dot{q},\ddot{q})=(\textrm{\scriptsize L}_1+\textrm{\scriptsize L}_3)\ddot{q}^3+\textrm{\scriptsize 
L}_1\ddot{q}^4-\f{q^1}{\textrm{\scriptsize C}_2}+\f{q^2}{\textrm{\scriptsize C}_2}+\f{q^3}{\textrm{\scriptsize 
C}_2}+const\nonumber\\
&&Q_4(q,\dot{q},\ddot{q})=\textrm{\scriptsize L}_1\ddot{q}^3+(\textrm{\scriptsize L}_1+\textrm{\scriptsize L}_4)\ddot{q}^4+\f 
{q^1}{\textrm{\scriptsize C}_1}+\f{q^4}{\textrm{\scriptsize
C}_1}+const\label{ex96} \ea The Birkhoffian  (\ref{ex96}) is
\textbf{conservative}. The function $E_{\omega}(q,\dot{q})$ is
given by (\ref{energie}), that is, \ba E_{\omega}(q,\dot{q})&=&\f
1{2}\textrm{\scriptsize L}_1(\dot{q}^3+\dot{q}^4)^2+\f
1{2}\textrm{\scriptsize L}_2(\dot{q}^2)^2+ \f
1{2}\textrm{\scriptsize L}_3(\dot{q}^3)^2+\f
1{2}\textrm{\scriptsize
L}_4(\dot{q}^4)^2+\frac{1}{2\textrm{\scriptsize C}_1}(q^1+
{q}^4)^2+\nonumber\\
&&
\frac{1}{2\textrm{\scriptsize C}_2}(-q^1+{q}^2+q^3)^2
+\f 1{2\textrm{\scriptsize C}_3}(q^1)^2+\sum^4_{j=1}(const)_jq^j
\label{ex97}
\ea

\noindent Because we are in  a situation where the  network has  one  loop which contains only
 capacitors, the Birkhoffian corresponding to 
(\ref{ex96}) is  \textbf{\textit{not} regular}. Indeed, the first
row of the matrix  $\left[\f{\pa Q_j}{\pa
\ddot{q}^i}\right]_{i,j=1,2,3,4}$  contains only zeros, therefore
$\textrm{det}\left[\f{\pa Q_j}{\pa
\ddot{q}^i}\right]_{i,j=1,2,3,4}=0$.\\
 As we have stated in
section 3, we can reduce the configuration space from dimension
$4$ to dimension $3$.
 Using the first equation in  (\ref{ex96}) we define
$\bar{M}_c\subset M_c$ by \be
\bar{M}_c=\{q=(q^1,q^2,q^3,q^4)\in M_c/\, \left(\frac{1}{\textrm{\scriptsize C}_1}+\f {1}{\textrm{\scriptsize C}_2}+\f 
{1}{\textrm{\scriptsize C}_3}\right)q^1-\f{q^2}{\textrm{\scriptsize C}_2}-\f{q^3}{\textrm{\scriptsize 
C}_2}+\f{q^4}{\textrm{\scriptsize C}_1}+const=0 \}
\ee
On the reduced configuration space $\bar{M}_c$, in the coordinate system  $\bar{q}^1:=q^2,\, \bar{q}^2:=q^3,\, 
\bar{q}^3:=q^4$, the Birkhoffian has the form
$\bar{\omega}_c=\sum^{3}_{j=1}\bar{Q}_jd\bar{q}^j$, \ba
&&\bar{Q}_1(\bar{q},\dot{\bar{q}},\ddot{\bar{q}})=\textrm{\scriptsize
L}_2\ddot{\bar{q}}^1+\mathfrak{C}_1
\bar{q}^1+\mathfrak{C}_1\bar{q}^2+\mathfrak{C}_2\bar{q}^3+const\nonumber\\
&&\bar{Q}_2(\bar{q},\dot{\bar{q}},\ddot{\bar{q}})=(\textrm{\scriptsize L}_1+\textrm{\scriptsize L}_3)\ddot{\bar{q}}^2+
\textrm{\scriptsize L}_1\ddot{\bar{q}}^3+\mathfrak{C}_1\bar{q}^1+
\mathfrak{C}_1\bar{q}^2+\mathfrak{C}_2\bar{q}^3+const\nonumber\\
&&\bar{Q}_3(\bar{q},\dot{\bar{q}},\ddot{\bar{q}})=\textrm{\scriptsize L}_1\ddot{\bar{q}}^2+
(\textrm{\scriptsize L}_1+\textrm{\scriptsize L}_4)\ddot{\bar{q}}^3+\mathfrak{C}_2\bar{q}^1
+\mathfrak{C}_2\bar{q}^2+\mathfrak{C}_3\bar{q}^3+const
\label{ex99}
\ea
where we have introduced the notation $\mathfrak{C}_1:=\f{1}{\textrm{\scriptsize C}_2}\left(1-\f{1}{\textrm{\scriptsize 
C}_2}\left(\frac{1}{\textrm{\scriptsize C}_1}+
\f {1}{\textrm{\scriptsize C}_2}+\f {1}{\textrm{\scriptsize C}_3}\right)^{-1}\right)$, \\
$\mathfrak{C}_2:=
\f{1}{\textrm{\scriptsize C}_2\textrm{\scriptsize C}_1}\left(\frac{1}{\textrm{\scriptsize C}_1}+\f {1}{\textrm{\scriptsize C}_2}+\f 
{1}{\textrm{\scriptsize C}_3}\right)^{-1}$,
 $\mathfrak{C}_3:=\f{1}{\textrm{\scriptsize C}_1}\left(1-\f{1}{\textrm{\scriptsize C}_1}\left(\frac{1}{\textrm{\scriptsize C}_1}+\f 
{1}{\textrm{\scriptsize C}_2}+
\f {1}{\textrm{\scriptsize C}_3}\right)^{-1}\right)$.

\noindent Let us now see whether the Birkhoffian (\ref{ex99}) is
\textbf{ regular} and/or  \textbf{conservative}. We calculate
 \be
\textrm{det}\left[\f{\pa\bar{Q}_j}{\pa \ddot{\bar{q}}^i}\right]_{i,j=1,2,3}=
\left|
\begin{array}{ccc}
\textrm{\scriptsize L}_2&0&0\\
0&\textrm{\scriptsize L}_1+\textrm{\scriptsize L}_3&\textrm{\scriptsize L}_1\\
0&\textrm{\scriptsize L}_1&\textrm{\scriptsize L}_1+\textrm{\scriptsize L}_4
\end{array}
\right|
\ee
Thus, if $\textrm{\scriptsize L}_2\left[\textrm{\scriptsize L}_1\textrm{\scriptsize L}_4+\textrm{\scriptsize 
L}_3(\textrm{\scriptsize L}_1+\textrm{\scriptsize L}_4)\right]\ne 0$, then the Birkhoffian (\ref{ex99}) is 
\textbf{regular}.

\noindent The corresponding Birkhoffian vector field (see section
2), is given by: \ba Y=&&\dot{\bar{q}}^1\f{\pa}{\pa
\bar{q}^1}+\dot{\bar{q}}^2\f{\pa}{\pa \bar{q}^2}+
\dot{\bar{q}}^3\f{\pa}{\pa \bar{q}^3} -\f{1}{\textrm{\scriptsize
L}_2}\left[\mathfrak{C}_1\bar{q}^1+
\mathfrak{C}_1\bar{q}^2+\mathfrak{C}_2\bar{q}^3\right]
\f{\pa}{\pa \dot{\bar{q}}^1}+\nonumber\\
&& \f {1}{(\textrm{\scriptsize L}_1+\textrm{\scriptsize L}_3)\textrm{\scriptsize L}_4+\textrm{\scriptsize L}_1\textrm{\scriptsize 
L}_3}
\left[\left(-\mathfrak{C}_1(\textrm{\scriptsize L}_1+\textrm{\scriptsize L}_4)+\mathfrak{C}_2\textrm{\scriptsize 
L}_1\right)\bar{q}^1+
\left(-\mathfrak{C}_1(\textrm{\scriptsize L}_1+\textrm{\scriptsize L}_4)+\mathfrak{C}_2\textrm{\scriptsize 
L}_1\right)\bar{q}^2+\right.\nonumber\\
&&\left.\left(-\mathfrak{C}_2(\textrm{\scriptsize L}_1+\textrm{\scriptsize L}_4)+\mathfrak{C}_3\textrm{\scriptsize 
L}_1\right)\bar{q}^3+const\right]
\f{\pa}{\pa \dot{q}^2}+
\nonumber\\
&& \f {1}{(\textrm{\scriptsize L}_1+\textrm{\scriptsize L}_4)\textrm{\scriptsize L}_3+\textrm{\scriptsize L}_1\textrm{\scriptsize 
L}_4}
\left[\left(-\mathfrak{C}_2(\textrm{\scriptsize L}_1+\textrm{\scriptsize L}_3)+\mathfrak{C}_1\textrm{\scriptsize 
L}_1\right)\bar{q}^1+
\left(-\mathfrak{C}_2(\textrm{\scriptsize L}_1+\textrm{\scriptsize L}_3)+\mathfrak{C}_1\textrm{\scriptsize 
L}_1\right)\bar{q}^2+\right.\nonumber\\
&&\left.\left(-\mathfrak{C}_3(\textrm{\scriptsize L}_1+\textrm{\scriptsize L}_3)+\mathfrak{C}_2\textrm{\scriptsize 
L}_1\right)\bar{q}^3+const\right] \f{\pa}{\pa \dot{q}^3} \ea Also,
the Birkhoffian  (\ref{ex99}) is \textbf{conservative} with the
 function $\bar{E}_{\bar{\omega}}(\bar{q},\dot{\bar{q}})$ given by \ba
\bar{E}_{\bar{\omega}}(\bar{q},\dot{\bar{q}})&=&\f 1{2}\textrm{\scriptsize L}_1(\dot{\bar{q}}^2+\dot{\bar{q}}^3)^2+\f 
1{2}\textrm{\scriptsize L}_2(\dot{\bar{q}}^1)^2+
\f 1{2}\textrm{\scriptsize L}_3(\dot{\bar{q}}^2)^2+\f 1{2}\textrm{\scriptsize 
L}_4(\dot{\bar{q}}^3)^2+\frac{1}{2}\mathfrak{C}_1(\bar{q}^1+
\bar{q}^2)^2+\nonumber\\
&&\mathfrak{C}_2(\bar{q}^1\bar{q}^3+\bar{q}^2\bar{q}^3)+\frac{1}{2}\mathfrak{C}_3(\bar{q}^3)^2
+\sum^3_{j=1}(const)_j\bar{q}^j \ea As we have pointed out in
section 3, because the network  has one   loop which contains only
inductors, this is the loop $I_3$, we can further reduce the
dimension of the configuration space by one. The equation that we
use for doing this is Kirchhoff's voltage law equation for this
loop, that is, the sixth equation in (\ref{ex92}). For the chosen
$q$-coordinate system  (\ref{csex9}) and after the transformation
$\mathcal{C}A^T$, the sixth equation in (\ref{ex92}) added with
the five equation in (\ref{ex92}) and it appeared in (\ref{ex96})
by the function $Q_3(q,\dot{q},\ddot{q})$. The Birkhoffian
formulation was presented in a coordinate free fashion.
In order to have a  coordinate system in which the sixth equation in (\ref{ex92}) appears 
in the initial form, we change the $\bar{q}$-coordinate system by
the following relations \ba
&&\bar{q}^1=\check{q}^1-\check{q}^2\nonumber\\
&& \bar{q}^2=\check{q}^2\nonumber\\
&& \bar{q}^3=\check{q}^3\label{coordonate^} \ea In terms of
$\check{q}$-coordinates on $\bar{M}_c$, the Birkhoffian
$\bar{\omega}_c=\sum^{3}_{j=1}\check{Q}_jd\check{q}^j$, where \ba
&&\check{Q}_1(\check{q},\dot{\check{q}},\ddot{\check{q}})=\textrm{\scriptsize L}_2\ddot{\check{q}}^1-\textrm{\scriptsize 
L}_2\ddot{\check{q}}^2+\mathfrak{C}_1\check{q}^1+\mathfrak{C}_2\check{q}^3+const\nonumber\\
&&\check{Q}_2(\check{q},\dot{\check{q}},\ddot{\check{q}})=-\textrm{\scriptsize L}_2\ddot{\check{q}}^1+(\textrm{\scriptsize 
L}_1+\textrm{\scriptsize L}_2+\textrm{\scriptsize
L}_3)\ddot{\check{q}}^2+
\textrm{\scriptsize L}_1\ddot{\check{q}}^3\nonumber\\
&&\check{Q}_3(\check{q},\dot{\check{q}},\ddot{\check{q}})=\textrm{\scriptsize
L}_1\ddot{\check{q}}^2+ (\textrm{\scriptsize
L}_1+\textrm{\scriptsize L}_4)\ddot{\check{q}}^3
+\mathfrak{C}_2\check{q}^1+\mathfrak{C}_3\check{q}^3+const
\label{ex910} \ea Using the second equation (\ref{ex910}), we
define $\hat{M}_{c}\subset \bar{M}_c$ by \be
\hat{M}_{c}=\{\check{q}=(\check{q}^1,\check{q}^2,\check{q}^3)\in
\bar{M}_c/\, -\textrm{\scriptsize
L}_2{\check{q}}^1+(\textrm{\scriptsize L}_1+\textrm{\scriptsize
L}_2+\textrm{\scriptsize L}_3){\check{q}}^2+ \textrm{\scriptsize
L}_1{\check{q}}^3+g(t)+const=0\} \ee with a certain function
$g(t)$ depending on t.

\noindent On the reduced configuration space $\hat{M}_{c}$, in the coordinate system $\hat{q}^1:=\check{q}^1,\, 
\hat{q}^2:=\check{q}^3$, the Birkhoffian has the form 
$\hat{\omega}=\hat{Q}_1d\hat{q}^1+\hat{Q}_2d\hat{q}^2$ with \ba
&&\hat{Q}_1(\hat{q},\dot{\hat{q}},\ddot{\hat{q}})=
\f{\textrm{\scriptsize L}_2(\textrm{\scriptsize L}_1+\textrm{\scriptsize L}_3)}{\textrm{\scriptsize L}_1+\textrm{\scriptsize 
L}_2+\textrm{\scriptsize L}_3}\ddot{\hat{q}}^1+\f{\textrm{\scriptsize L}_1\textrm{\scriptsize L}_2}{\textrm{\scriptsize 
L}_1+\textrm{\scriptsize L}_2+\textrm{\scriptsize 
L}_3}\ddot{\hat{q}}^2+\mathfrak{C}_1\hat{q}^1+\mathfrak{C}_2\hat{q}^2+const
\nonumber\\
&&\hat{Q}_2(\hat{q},\dot{\hat{q}},\ddot{\hat{q}})=
\f{\textrm{\scriptsize L}_1\textrm{\scriptsize L}_2}{\textrm{\scriptsize L}_1+\textrm{\scriptsize L}_2+\textrm{\scriptsize 
L}_3}\ddot{\hat{q}}^1+
\f{(\textrm{\scriptsize L}_1+\textrm{\scriptsize L}_4)(\textrm{\scriptsize L}_2+\textrm{\scriptsize L}_3)+\textrm{\scriptsize 
L}_1\textrm{\scriptsize L}_4}{\textrm{\scriptsize
L}_1+\textrm{\scriptsize L}_2+\textrm{\scriptsize
L}_3}\ddot{\hat{q}}^2
+\mathfrak{C}_2\hat{q}^1+\mathfrak{C}_3\hat{q}^2+const
\nonumber\\
&&\label{maisus}
\ea
Because ${\textrm{\scriptsize L}_i}\neq0$, $i=0,...,4$,   the determinant $\f{\textrm{\scriptsize L}_2}{(\textrm{\scriptsize L}_1+\textrm{\scriptsize L}_2+\textrm{\scriptsize 
L}_3)^2}\left|
\begin{array}{cc}
\textrm{\scriptsize L}_1+\textrm{\scriptsize L}_3&\textrm{\scriptsize L}_1\\
\textrm{\scriptsize L}_1\textrm{\scriptsize L}_2&
(\textrm{\scriptsize L}_1+\textrm{\scriptsize L}_4)(\textrm{\scriptsize L}_2+\textrm{\scriptsize L}_3)+\textrm{\scriptsize 
L}_1\textrm{\scriptsize L}_4
\end{array}
\right|\ne 0$, then, the Birkhoffian (\ref{maisus}) is
\textbf{regular}.

\noindent Moreover, the Birkhoffian (\ref{maisus}) is
\textbf{conservative}, and \ba
\hat{E}_{\omega}(\hat{q},\dot{\hat{q}})&=&
\f{\textrm{\scriptsize L}_2(\textrm{\scriptsize L}_1+\textrm{\scriptsize L}_3)}{2(\textrm{\scriptsize L}_1+\textrm{\scriptsize 
L}_2+\textrm{\scriptsize L}_3)}(\dot{\hat{q}}^1)^2+\f{\textrm{\scriptsize L}_1\textrm{\scriptsize L}_2}{\textrm{\scriptsize 
L}_1+\textrm{\scriptsize L}_2+\textrm{\scriptsize
L}_3}\dot{\hat{q}}^1\dot{\hat{q}}^2+
\f{(\textrm{\scriptsize L}_1+\textrm{\scriptsize L}_4)(\textrm{\scriptsize L}_2+\textrm{\scriptsize L}_3)+\textrm{\scriptsize 
L}_1\textrm{\scriptsize L}_4}{2(\textrm{\scriptsize L}_1+\textrm{\scriptsize L}_2+\textrm{\scriptsize 
L}_3)}(\dot{\hat{q}}^2)^2+
+\nonumber\\
&&\f{1}{2}\mathfrak{C}_1(\hat{q}^1)^2+\mathfrak{C}_2\hat{q}^1\hat{q}^2+\frac{1}{2}\mathfrak{C}_3(\hat{q}^2)^
2 +\sum^2_{j=1}(const)_j\hat{q}^j \ea

Let us now suppose that the inductors and the capacitors in the network  are nonlinear devices,  their constitutive relations
being of the form (\ref{2}), 
(\ref{4}).
The equations (\ref{7}) which govern  the network have now the form
\be
\left\{\begin{array}{lllllllllll}
\textsc{i}_4-\f{d\textsc{q}_1}{dt}+\f{d\textsc{q}_3}{dt}=0\\
\textsc{i}_2+\textsc{i}_3-\f{d\textsc{q}_2}{dt}-\f{d\textsc{q}_3}{dt}=0\\
\textsc{i}_1-\textsc{i}_3-\textsc{i}_4=0\\
\\
C_1(\textsc{q}_1)-C_2(\textsc{q}_2)+C_3(\textsc{q}_3)=0\\
L_2(\textsc{i}_2)\f{d\textsc{i}_2}{dt}+
C_2(\textsc{q}_2)=0\\
L_1(\textsc{i}_1)\f{d\textsc{i}_1}{dt}-
L_2(\textsc{i}_2)\f{d\textsc{i}_2}{dt}+
L_3(\textsc{i}_3)\f{d\textsc{i}_3}{dt}=0\\
-L_1(\textsc{i}_1)\f{d\textsc{i}_1}{dt}-
L_4(\textsc{i}_4)\f{d\textsc{i}_4}{dt}-C_1(\textsc{q}_1)=0
\end{array}
\right.\label{exn92}
\ee
where  $L_a:\mathbf{R}\longrightarrow \mathbf{R}\backslash \{0\}$,  $C_\alpha:\mathbf{R}\longrightarrow \mathbf{R}\backslash\{0\}$ 
are smooth invertible functions.

\noindent As we have pointed out in section 3,  the first set of equations in  (\ref{exn92}) is the same as in the linear case, therefore the 
configuration  space $M_c$ is the same, too.
For the  coordinate system on  $M_c$  given by (\ref{csex9}),
the matrices  $\mathcal{N}$, $\mathcal{C}$  have the same expressions (\ref{N9}) as before. Thus,  in  the nonlinear case, the  
Birkhoffian becomes  $\omega_c=\sum^{4}_{j=1}Q_j(q,\dot{q},\ddot{q})dq^j$  with the functions $Q_j$ given by (\ref{bir}), 
(\ref{bir2}), (\ref{bir1}), that is, \ba
Q_1(q,\dot{q},\ddot{q})&=&C_3(q^1)+C_1(q^1+q^4+const)-C_2(-q^1+q^2+q^3+const)\nonumber\\
Q_2(q,\dot{q},\ddot{q})&=&L_2(\dot{q}^2)\ddot{q}^2+C_2(-q^1+q^2+q^3+const)\nonumber\\
Q_3(q,\dot{q},\ddot{q})&=&\left(L_3(\dot{q}^3)+
L_1(\dot{q}^3+\dot{q}^4)\right)\ddot{q}^3+
L_1(\dot{q}^3+\dot{q}^4)\ddot{q}^4+\nonumber\\
&&C_2(-q^1+q^2+q^3+const)
\nonumber\\
Q_4(q,\dot{q},\ddot{q})&=&
L_1(\dot{q}^3+\dot{q}^4)\ddot{q}^3+\left(L_4(\dot{q}^4)+L_1(\dot{q}^3+\dot{q}^4)\right)\ddot{q}^4+
C_1(q^1+q^4+const)\nonumber\\
\label{exn96} \ea The Birkhoffian  (\ref{exn96}) is
\textbf{conservative} with  the  function $E_{\omega}(q,\dot{q})$
given by
 (\ref{energienonlinear}), that is,
\ba E_{\omega}(q,\dot{q})&=& \int
\widetilde{L}_1(\dot{q})(\dot{q}^3+\dot{q}^4)(d\dot{q}^3+d\dot{q}^4)+
\int L_2(\dot{q}^2)\dot{q}^2d\dot{q}^2+
\int L_3(\dot{q}^3)\dot{q}^3d\dot{q}^3+\nonumber\\
&& \int L_4(\dot{q}^4)\dot{q}^4d\dot{q}^4- \int \int
\widetilde{L}_1'(\dot{q})(\dot{q}^3+\dot{q}^4)d\dot{q}^3d\dot{q}^4-
\int \int \widetilde{L}_1(\dot{q})d\dot{q}^3d\dot{q}^4+\nonumber\\
&& \int \widetilde{C}_1(q)(dq^1+dq^4)+ \int
\widetilde{C}_2(q)(-dq^1+dq^2+dq^3)+
\int C_3(q^1)dq^1-\nonumber\\
&& \int \int \widetilde{C}_1'(q)dq^1dq^4- \int \int
\widetilde{C}_2'(q)(-dq^1dq^2-dq^1dq^3+dq^2dq^3)
\nonumber\\
&&-\int\int\int \widetilde{C}_2''(q)dq^1dq^2dq^3 \label{exn97} \ea
The Birkhoffian  (\ref{exn96}) is \textbf {\textit{not}  regular}, since the first row of the matrix $\left[\f{\pa Q_j}{\pa 
\ddot{q}^i}\right]_{i,j=1,2,3,4}$ contains only zeros.
But just as in the linear case, we can reduce the configuration space from dimension $4$ to dimension $3$. Different from the linear 
case, the reduced configuration space will {\it \textbf{not}} be a linear subspace of $M_c$.

\noindent If  the functions $C_1$, $C_2$, $C_3$ are such that the
Jacobian matrix for the first equation in  (\ref{exn96})  has rank
one, we define the 3-dimensional manifold $\bar{M}_c\subset M_c$
by \be \bar{M}_c=\{q=(q^1,q^2,q^3,q^4)\in M_c/\,
C_3(q^1)+C_1(q^1+q^4+const)-C_2(-q^1+q^2+q^3+const)=0
\}\label{barM} \ee
By the implicit function theorem, we obtain a  local coordinate system on the reduced configuration space $\bar{M}_c$. Taking 
$\bar{q}^1:=q^2,\, \bar{q}^2:=q^3,\, \bar{q}^3:=q^4$, the Birkhoffian has the form 
$\bar{\omega}_c=\sum^{3}_{j=1}\bar{Q}_jd\bar{q}^j$, with
\ba
\bar{Q}_1(\bar{q},\dot{\bar{q}},\ddot{\bar{q}})&=&
L_2(\dot{\bar{q}}^1)\ddot{\bar{q}}^1+C_2(-f(\bar{q}^1,\bar{q}^2,\bar{q}^3)+\bar{q}^1+\bar{q}^2+const)\nonumber\\
\bar{Q}_2(\bar{q},\dot{\bar{q}},\ddot{\bar{q}})&=&
\left(L_3(\dot{\bar{q}}^2)+
L_1(\dot{\bar{q}}^2+\dot{\bar{q}}^3)\right)\ddot{\bar{q}}^2+
L_1(\dot{\bar{q}}^2+\dot{\bar{q}}^3)\ddot{\bar{q}}^3+\nonumber\\
&&C_2(-f(\bar{q}^1,\bar{q}^2,\bar{q}^3)+\bar{q}^1+\bar{q}^2+const)\nonumber\\
\bar{Q}_3(\bar{q},\dot{\bar{q}},\ddot{\bar{q}})&=&L_1(\dot{\bar{q}}^2+\dot{\bar{q}}^3)\ddot{\bar{q}}^2+
\left(L_4(\dot{\bar{q}}^3)+
L_1(\dot{\bar{q}}^2+\dot{\bar{q}}^3)\right)
\ddot{\bar{q}}^3+\nonumber\\
&&C_1(f(\bar{q}^1,\bar{q}^2,\bar{q}^3)+\bar{q}^3+const)
\label{exn99}
\ea
 $f:U\subset\mathbf{R}^3\longrightarrow \mathbf{R}^1$ being an unique function such that
$f(\bar{q}_0)=q^1_0$, $q^1_0\in \mathbf{R}$, and
$C_3(f(\bar{q}))+C_1(f(\bar{q})+\bar{q}^3+const)-C_2(-f(\bar{q})+\bar{q}^1+\bar{q}^2+const)=0$, $\forall 
\bar{q}=(\bar{q}^1,\bar{q}^2,\bar{q}^3)\in U$, with $U$ a
neighborhood of $\bar{q}_0=(\bar{q}_0^1,\bar{q}_0^2,\bar{q}_0^3)$.
On account of
$L_1,L_2,L_3,L_4:\mathbf{R}\longrightarrow\mathbf{R}\backslash
\{0\}$, we have \be \left|
\begin{array}{ccc}
L_2(\dot{\bar{q}}^1)&0&0\\
0&L_3(\dot{\bar{q}}^2)+
L_1(\dot{\bar{q}}^2+\dot{\bar{q}}^3)&L_1(\dot{\bar{q}}^2+\dot{\bar{q}}^3)\\
0&L_1(\dot{\bar{q}}^2+\dot{\bar{q}}^3)&L_4(\dot{\bar{q}}^3)+
L_1(\dot{\bar{q}}^2+\dot{\bar{q}}^3)
\end{array}
\right|\ne 0 \ee then, the Birkhoffian  (\ref{exn99}) is
\textbf{regular}.

\noindent Because the network  has one   loop which contains only
inductors, let us perform a further reduction of  the dimension of
the configuration space by one, just as we have done in the linear
case. In the coordinate system $\check{q}$ defined in
(\ref{coordonate^}), the Birkoffian
$\bar{\omega}_c=\sum^{3}_{j=1}\check{Q}_jd\check{q}^j$, where \ba
\check{Q}_1(\check{q},\dot{\check{q}},\ddot{\check{q}})&=&
L_2(\dot{\check{q}}^1-\dot{\check{q}}^2)\ddot{\check{q}}^1-
L_2(\dot{\check{q}}^1-\dot{\check{q}}^2)
\ddot{\check{q}}^2+C_2(\check{q}^1,\check{q}^2,\check{q}^3)\nonumber\\
\check{Q}_2(\check{q},\dot{\check{q}},\ddot{\check{q}})&=&-
L_2(\dot{\check{q}}^1-\dot{\check{q}}^2)\ddot{\check{q}}^1+
\left[L_1(\dot{\check{q}}^2+\dot{\check{q}}^3)+
L_2(\dot{\check{q}}^1-\dot{\check{q}}^2)+
L_3(\dot{\check{q}}^2)\right]\ddot{\check{q}}^2+\nonumber\\
&&
L_1(\dot{\check{q}}^2+\dot{\check{q}}^3)\ddot{\check{q}}^3\nonumber\\
\check{Q}_3(\check{q},\dot{\check{q}},\ddot{\check{q}})&=&L_1(\dot{\check{q}}^2+\dot{\check{q}}^3)\ddot{\check{q}}^2+
\left[L_1(\dot{\check{q}}^2+\dot{\check{q}}^3)+
L_4(\dot{\check{q}}^3)\right]\ddot{\check{q}}^3
+C_1(\check{q}^1,\check{q}^2,\check{q}^3)\nonumber\\
&&
\label{exn910}
\ea
Using the second equation in (\ref{exn910}),
 we can define a smooth constant rank affine sub-bundle $\mathfrak{S}_c$ of the affine bundle $\pi_J:J^2(\bar{M}_c)\longrightarrow 
T\bar{M}_c$ via \ba \mathfrak{S}_c&=&\{(\check{q},
\dot{\check{q}}, \ddot{\check{q}}) \in J^2(\bar{M_c})/\, -
L_2(\dot{\check{q}}^1-\dot{\check{q}}^2)\ddot{\check{q}}^1+
\left[L_1(\dot{\check{q}}^2+\dot{\check{q}}^3)+
L_2(\dot{\check{q}}^1-\dot{\check{q}}^2)+
L_3(\dot{\check{q}}^2)\right]\ddot{\check{q}}^2+\nonumber\\
&& L_1(\dot{\check{q}}^2+\dot{\check{q}}^3)\ddot{\check{q}}^3 =0\}
\ea The constraint $\mathfrak{S}_c$ is integrable, in the sense
that we have the foliation \be \mathfrak{F}_{const}:=\{(\check{q},
\dot{\check{q}})\in T(\bar{M_c})/\,
-\mathcal{L}_2(\dot{\check{q}}^1-\dot{\check{q}}^2)+\mathcal{L}_1(\dot{\check{q}}^2+\dot{\check{q}}^3)+\mathcal{L}_3(\dot{\check{q}}^2)=const\}\label{exn911}
\ee Thus, in the nonlinear case, we draw the conclusion that  we
can further reduce the configuration space
only if it is possible to 
find from (\ref{exn911}) new configuration coordinates $\hat{q}^1$, $\hat{q}^2$, that is, when the constraint 
(\ref{exn911}) is holonomic.

\vspace{0.5cm}

In order to underline that, depending on the topology of the networks with independent sources, the associated Birkhoffian is 
conservative or not,  we consider the circuit shown in figure 2 below. This circuit  contains  a  loop formed by capacitors and independent 
voltage sources and a  cutset formed by inductors and independent current sources. We shall see that the Birkhoffian associated to 
such a circuit is \textbf{\textit{not} regular} and  \textbf{\textit{not} {even} conservative} for nonlinear inductors and 
capacitors.

\begin{center}
\scalebox{0.40} {\includegraphics*{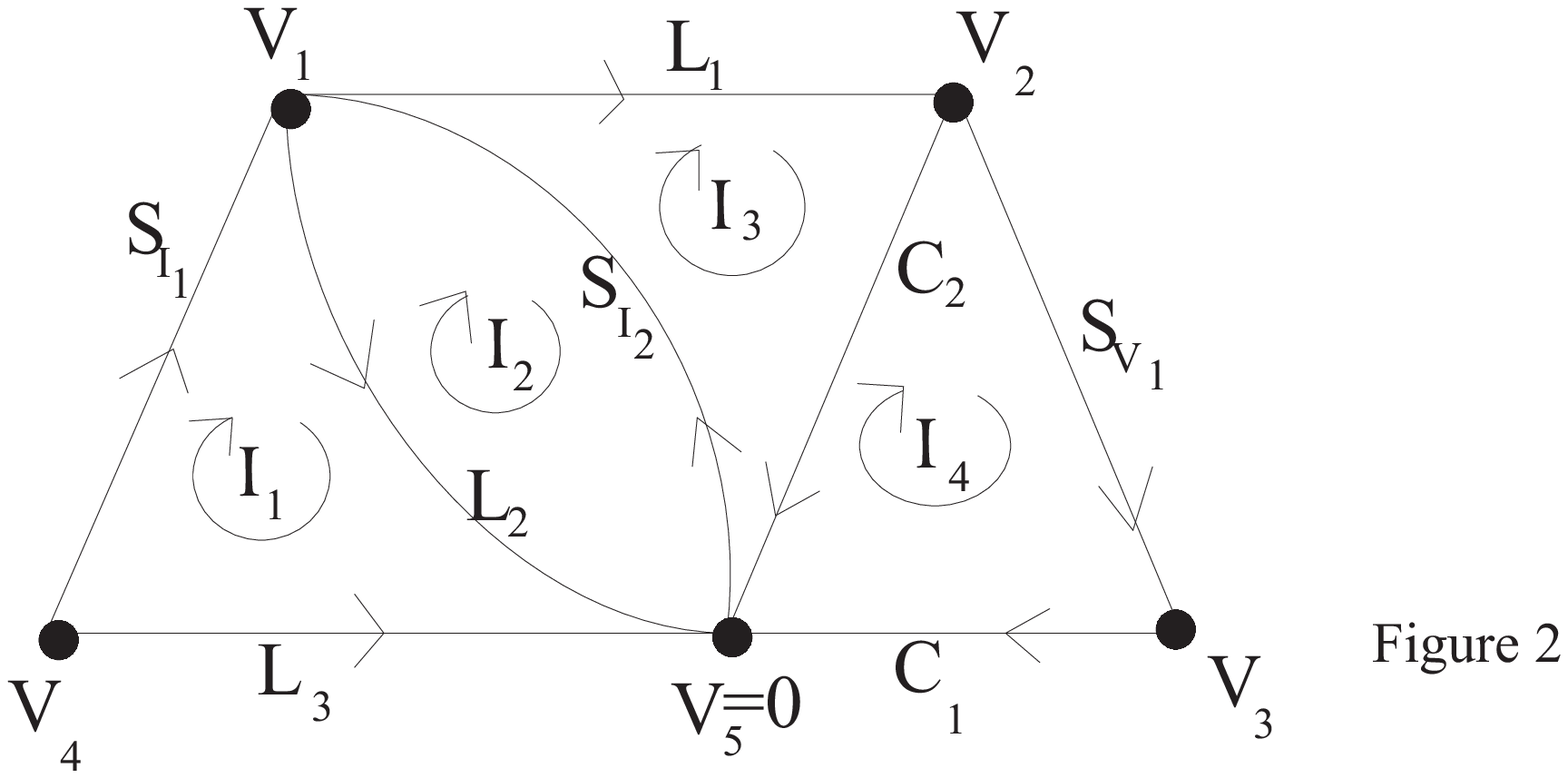}}
\end{center}

\noindent We have $k=3$, $p=2$, ${\textrm {\scriptsize S}_I}=2$,
${\textrm {\scriptsize S}_V}=1$, $n=4$, $m=4$ , $b=8$. We choose
the reference node to be $V_5$ and the current directions as
indicated in Figure 2. We cover the associated graph with the
loops
 $I_1$, $I_2$, $I_3$, $I_4$. Let $V=(V_1,V_2,V_3,V_4)\in \mathbf{R}^4$ be the vector of node voltage values,
$\textsc{i}=(\textsc{i}_{a},\textsc{i}_\alpha,\textsc{i}_{\textrm
{\scriptsize S}_I}, \textsc{i}_{\textrm {\scriptsize S}_V})\in
\mathbf{R}^3\times \mathbf{R}^2\times \mathbf{R}^2\times
\mathbf{R}^1$ be
 the vector of branch current  values  and
 $v=(v_{a},v_{\alpha},v_{\textrm {\scriptsize S}_I},{\textrm {\scriptsize S}_V})\in
\mathbf{R}^3\times \mathbf{R}^2\times \mathbf{R}^2\times
\mathbf{R}^1 $ be
 the vector of branch  voltage
 values.\\
The branches in Figure 2 are  labelled as follows: the first, the
second, and the third  branch are the inductor branches
$\textsc{L}_1$, $\textsc{L}_2$, $\textsc{L}_3$, the forth and the
fifth branch are the capacitor branches  $\textsc{C}_1$,
$\textsc{C}_2$, the next two branches are the current source
branches $\textsc{S}_{\textsc{I}_1}$, $\textsc{S}_{\textsc{I}_2}$,
and the last branch is the voltage source branch
$\textsc{S}_{\textsc{V}_1}$. The incidence and loop matrices,
$B\in \mathfrak{M}_{84}(\mathbf{R})$ and
 $A\in \mathfrak{M}_{84}(\mathbf{R})$, write as
\be  B=\left(
\begin{array}{cccc}
-1& 1& 0& 0\\
-1& 0& 0& 0\\
 0& 0& 0&-1 \\
 0& 0&-1& 0\\
0&-1& 0& 0\\
1& 0& 0&-1 \\
1& 0& 0& 0\\
0&-1& 1& 0
\end{array}
\right),\quad A=\left(
\begin{array}{cccc}
 0& 0& 1& 0 \\
 1& -1& 0& 0\\
-1& 0& 0& 0\\
 0& 0& 0& 1 \\
 0& 0& 1&-1\\
 1& 0& 0& 0 \\
 0& -1& 1& 0\\
 0& 0& 0& 1
\end{array}
\right)
 \label{exs21}
\ee For linear inductors and capacitors, the governing equations
have the form: \be \left\{
\begin{array}{lllllllll}
-\textsc{i}_1-\textsc{i}_2+\textsc{i}_{s_{I_1}}(t)+\textsc{i}_{s_{I_2}}(t)=0\\
\textsc{i}_1-\f{d\textsc{q}_2}{dt}-\textsc{i}_{s_{V_1}}=0\\
-\f{d\textsc{q}_1}{dt}+\textsc{i}_{s_{V_1}}=0\\
-\textsc{i}_3-\textsc{i}_{s_{I_1}}(t)=0\\
\\
{\textrm{\scriptsize{L}}}_2\f{d\textsc{i}_2}{dt}-{\textrm{\scriptsize{L}}}_3\f{d\textsc{i}_3}{dt}+v_{s_{I_1}}=0\\
-{\textrm{\scriptsize{L}}}_2\f{d\textsc{i}_2}{dt}-v_{s_{I_2}}=0\\
{\textrm{\scriptsize{L}}}_1\f{d\textsc{i}_1}{dt}+\f
{\textsc{q}_2}{{\textrm{\scriptsize{C}}}_2}+
v_{s_{I_2}}=0\\
\f {\textsc{q}_1}{{\textrm{\scriptsize{C}}}_1}-\f {\textsc{q}_2}{{\textrm{\scriptsize{C}}}_2}+v_{s_{V_1}}(t)=0
\end{array}
\right.\label{exs22}
\ee
where ${\textrm{\scriptsize{C}}}_\alpha\neq 0$, $\alpha=1,2$ and ${\textrm{\scriptsize{L}}}_a\neq 0$, $a=1,2,3$, are distinct 
constants.

\noindent Note that  $\textsc{i}_{s_{I_1}}$, $\textsc{i}_{s_{I_2}}$,  $v_{s_{V_1}}$ are given functions of time which describe the currents 
associated to the independent current sources $S_{I_1}$, $S_{I_2}$ and the voltage associated to the independent voltage source 
$S_{V_1}$,  respectively.

\noindent Once we know the unknowns  $\textsc{i}_1$, $\textsc{i}_2$, $\textsc{i}_3$, $\textsc{q}_1$, $\textsc{q}_2$, 
we can determine all the other circuit variables.

\noindent From the first set of equations $(\ref{exs22})$, we have
\be \textsc{i}_{s_{V_1}}=\f{d\textsc{q}_1}{dt} \ee and from the
second set of equations (\ref{exs22}), we conclude \ba
v_{s_{I_1}}=-{\textrm{\scriptsize{L}}}_2\f{d\textsc{i}_2}{dt}+{\textrm{\scriptsize{L}}}_3\f{d\textsc{i}_3}{dt}\nonumber\\
v_{s_{I_2}}=-{\textrm{\scriptsize{L}}}_2\f{d\textsc{i}_2}{dt} \ea
Therefore,  the system (\ref{surse7*}) has now the form \be
\left\{
\begin{array}{lllllllll}
-\textsc{i}_1-\textsc{i}_2+\textsc{i}_{s_{I_1}}(t)+\textsc{i}_{s_{I_2}}(t)=0\\
\textsc{i}_1-\f{d\textsc{q}_2}{dt}-\f{d\textsc{q}_1}{dt}=0\\
-\textsc{i}_3-\textsc{i}_{s_{I_1}}(t)=0\\
\\
{\textrm{\scriptsize{L}}}_1\f{d\textsc{i}_1}{dt}+\f 
{\textsc{q}_2}{{\textrm{\scriptsize{C}}}_2}-{\textrm{\scriptsize{L}}}_2\f{d\textsc{i}_2}{dt}=0\\
\f {\textsc{q}_1}{{\textrm{\scriptsize{C}}}_1}-\f {\textsc{q}_2}{{\textrm{\scriptsize{C}}}_2}+v_{s_{V_1}}(t)=0
\end{array}
\right.\label{exs22*}
\ee
with
\be
\mathcal{B}_1^T=
\left(
\begin{array}{ccccc}
-1& -1&0& 0&0 \\
 1& 0& 0&-1&-1\\
 0& 0& -1& 0&0
\end{array}
\right)\quad
\mathcal{B}_2^T=
\left(
\begin{array}{cc}
 1&1 \\
 0&0\\
-1&0
\end{array}
\right)
\ee
\be
\mathcal{A}_1^T=
\left(
\begin{array}{ccccc}
 1& -1&0&0& 1 \\
 0& 0& 0&1&- 1
\end{array}
\right)\quad
\mathcal{A}_2^T=
\left(
\begin{array}{c}
 0 \\
 1
\end{array}
\right) \ee The relations  (\ref{var}), (\ref{cssurse}), read as
follows for this example \be
\textsc{i}_a:=\f{d\textsc{q}_{(a)}}{dt}, \, a=1,2,3 \label{vars21}
\ee \be
x^1:=\textsc{q}_{(1)},\,x^2:=\textsc{q}_{(2)},\,x^{3}:=\textsc{q}_{(3)},\,
,x^{4}:=\textsc{q}_{1},\,x^5:=\textsc{q}_{2}\label{css21} \ee
Using the first set of equations (\ref{exs22*}), we  define the
2-dimensional affine-linear configuration  space $M_c$. We solve
the corresponding equations (\ref{surse7''})  in terms of 2
variables. In view of  the notations (\ref{vars21}),
(\ref{css21}), we obtain, for example, \ba
x^{2}&=&-x^1+f^2(t)+const\nonumber\\
x^{3}&=&f^3(t)+const\nonumber\\
x^5&=&x^1-x^4+const
\label{exs23}
\ea
with $f^2(t)=\int\left(\textsc{i}_{s_{1}}(t)+
\textsc{i}_{s_{2}}(t)\right){\textrm dt}$, $f^3(t)=-\int\textsc{i}_{s_{1}}(t){\textrm dt}$ and
the other components of $f$  in  (\ref{surse10'}) being zero.
Thus a coordinate system on $M_c$ is given by
\be
q^1:=x^1,q^2:=x^4.\label{csexs1}
\ee
and the matrices of constants $\mathfrak{N}=\left(\begin{array}{c}
\mathfrak{N}^{a}_{ j}\\
\mathfrak{N}^{\alpha}_{j}
\end{array}\right)_{{a=1,2,\, \alpha=3,4,5 \atop j=1,2}}$ and $\mathcal{C}$  are
\be
\mathfrak{N}=
\left(
\begin{array}{cc}
 1& 0 \\
 -1& 0\\
0&0\\
0&1\\
1&-1
\end{array}\right)
,\quad \mathcal{C}=\left(
\begin{array}{cc}
 1& 0\\
 0& 1
\end{array}
\right) \ee In terms of the coordinates (\ref{csexs1}), we define
the Birkhoffian  $\omega_{t_c}=Q_1dq^1+Q_2dq^2$,  as in
(\ref{sursebir})-(\ref{surse14}), that is, \ba
&&Q_1(t,q,\dot{q},\ddot{q})=({\textrm{\scriptsize{L}}}_1+{\textrm{\scriptsize{L}}}_2)\ddot{q}^1-{\textrm{\scriptsize{L}}}_2\f{d^2f^2
(t)}{dt^2}+\f{q^1}{{\textrm{\scriptsize{C}}}_2}-\f{q^2}{{\textrm{\scriptsize{C}}}_2}+const\nonumber\\
&&Q_2(t,q,\dot{q},\ddot{q})=-\f {q^1}{{\textrm{\scriptsize{C}}}_2}
+\left(\f{1}{{\textrm{\scriptsize{C}}}_1}+\f 1{{\textrm{\scriptsize{C}}}_2}\right)q^2+v_{s_{V_1}}(t)+const\nonumber\\
\label{exs27} \ea Because there exists a  loop which contains only
capacitors and independent voltage sources, the Birkhoffian
(\ref{exs27}) is   \textbf{\textit {not} regular}. Indeed, the
second row of the matrix $\left[\f{\pa Q_j}{\pa
\ddot{q}^i}\right]_{i,j=1,2}$ contains only zeros, therefore,
$\textrm{det}\left[\f{\pa Q_j}{\pa
\ddot{q}^i}\right]_{i,j=1,2}=0$.

Though there exists in the network a cutset formed by inductors and independent current sources,  the Birkhoffian 
(\ref{exs27}) is
\textbf{ conservative} in the sense of
 definition (\ref{surseconserv}).
The   function $E_{\omega_{t}}(t,q,\dot{q})$ is given by
(\ref{energiesurse1}), that is, \ba E_{\omega_{t}}(t,q,\dot{q})&=&
\f 1{2}\textrm{\scriptsize L}_1(\dot{q}^1)^2+\f
1{2}\textrm{\scriptsize L}_2(\dot{q}^1)^2+
\frac{1}{2\textrm{\scriptsize C}_1}(q^2)^2+
\frac{1}{2\textrm{\scriptsize C}_2}(q^1-q^2)^2+\nonumber\\
&&+\left(-{\textrm{\scriptsize{L}}}_2\f{d^2f^2(t)}{dt^2}+const_1\right)q^1+\left(v_{s_{V_1}}(t)+const_2\right)q^2
\ea

\noindent In order to obtain a  regular  Birkhoffian we could use
the second equation from (\ref{exs27}) and reduce the
configuration space $M_c$ to a vector space $\bar{M}_c$ of
dimension 1. The procedure is the same as in the first example
with linear devices.

For nonlinear inductors and capacitors, in the coordinate system
(\ref{csexs1}) on the configuration space $M_c$ of dimension $2$,
the Birkhoffian $\omega_{t_c}=Q_1dq^1+Q_2dq^2$, where
 \ba
Q_1(t,q,\dot{q},\ddot{q})&=&\left[L_1(\dot{q}^1)+L_2\left(-\dot{q}^1+\f{df^2(t)}{dt}
\right)\right]\ddot{q}^1-L_2\left(-\dot{q}^1+\f{df^2(t)}{dt}\right)\f{d^2f^2(t)}{dt^2}
+\nonumber\\
&&C_2(q^1-q^2+const)
\nonumber\\
Q_2(t,q,\dot{q},\ddot{q})&=&C_1(q^2)-C_2(q^1-q^2+const)+
v_{s_{V_1}}(t) \label{exsn26} \ea

\noindent The Birkhoffian  (\ref{exsn26}) is \textbf{ \textit{not}
regular} and \textbf{ \textit{not} conservative}. Indeed, two of
the necessary conditions for the existence of the function
$E_{\omega_t}:TM\to \mathbf{R}$ such that
 $\sum^{2}_{j=1}Q_j(t,q,\dot{q},\ddot{q})dq^j=
\sum^{2}_{j=1}\f{\pa E_{\omega_t}}{\pa q^j}\dot{q}^j+
\f{\pa E_{\omega_t}}{\pa \dot{q}^j}\ddot{q}^j$, are
\be
\left\{
\begin{array}{lll}
\f{\pa E_{\omega_t}}{\pa \dot{q}^1}=L_1(\dot{q}^1)+L_2\left(-\dot{q}^1+\f{df^2(t)}{dt^2}
\right)\\
\\
\f{\pa E_{\omega_t}}{\pa q^1}=-L_2\left(-\dot{q}^1+\f{df^2(t)}{dt^2}\right)\f{d^2f^2(t)}{dt^2}+C_2(q^1-q^2+const)
\end{array}\right.
\ee
We easily see that for almost all values of the parameters, $\f{\partial ^2 E_{\omega_t}}{\partial \dot{q}^1q^1}\ne \f{\partial^2 
E_{\omega_t}}{\partial q^1\dot{q}^1}=0$.
\\

Let us  now consider a network  that has the oriented connected
graph as in  Figure 2 in which we interchanged  the inductor
branch $\textsc{L}_3$ with the capacitor branch $\textsc{C}_1$ and
the inductor branch $\textsc{L}_2$  with the capacitor branch
$\textsc{C}_2$. We will see that the Birkhoffian associated  to
this network is \textbf{\textit{conservative}} even if the
inductors and the capacitors in the network are nonlinear devices.

\noindent The system (\ref{surse7*}) has now the form \be
\left\{
\begin{array}{lllllllll}
-\textsc{i}_1-\f{d\textsc{q}_2}{dt}+\textsc{i}_{s_{I_1}}(t)+\textsc{i}_{s_{I_2}}(t)=0\\
\textsc{i}_1-\textsc{i}_2-\textsc{i}_3=0\\
-\f{d\textsc{q}_1}{dt}-\textsc{i}_{s_{I_1}}(t)=0\\
\\
L_1(\textsc{i}_1)\f{d\textsc{i}_1}{dt}+L_2(\textsc{i}_2)\f{d\textsc{i}_2}{dt}-C_2(\textsc{q}_2)=0\\
-L_2(\textsc{i}_2)\f{d\textsc{i}_2}{dt}+L_3(\textsc{i}_3)\f{d\textsc{i}_3}{dt}+v_{s_{V_1}}(t)=0
\end{array}
\right.
\ee
and \be \textsc{i}_{s_{V_1}}=\textsc{i}_3\ee  \ba
v_{s_{I_1}}=C_1(\textsc{q}_1)-C_2(\textsc{q}_2)\nonumber\\
v_{s_{I_2}}=-C_2(\textsc{q}_2) \ea

\noindent Using the same procedure  as above we get  the
configuration space $\mathcal{M}_c$ of dimension $2$. In view of
the notations  (\ref{vars21}), (\ref{css21}), a coordinate system
on $\mathcal{M}_c$ is given by \be
q^1:=x^1,q^2:=x^3.\label{csexs1'} \ee and the Birkhoffian
$\mathcal{\omega}_{t_c}=Q_1dq^1+Q_2dq^2$, where
 \ba
Q_1(t,q,\dot{q},\ddot{q})&=&\left[L_1(\dot{q}^1)+L_2\left(\dot{q}^1-\dot{q}^2
\right)\right]\ddot{q}^1-L_2\left(\dot{q}^1-\dot{q}^2\right)\ddot{q}^2
-C_2(-q^1+f^5(t)+const)
\nonumber\\
Q_2(t,q,\dot{q},\ddot{q})&=&-L_2\left(\dot{q}^1-\dot{q}^2
\right)\ddot{q}^1+\left[L_2\left(\dot{q}^1-\dot{q}^2\right)+L_3(\dot{q}^2)\ddot{q}^2\right]\ddot{q}^2
+ v_{s_{V_1}}(t) \label{exsn26'} \ea with
$f^5(t)=\int\left(\textsc{i}_{s_{1}}(t)+
\textsc{i}_{s_{2}}(t)\right){\textrm dt}$ and the other components
of $f$  in  (\ref{surse10'}) being zero. The Birkhoffian
(\ref{exsn26'}) is conservative in the sense of
 definition (\ref{surseconserv}), with the  function
$E_{\omega_{t}}(t,q,\dot{q})$  given by  \ba
E_{\omega_{t}}(t,q,\dot{q})&=& \int
L_1(\dot{q}^1)\dot{q}^1d\dot{q}^1+ \int
\widetilde{L}_2(\dot{q})(\dot{q}^1-\dot{q}^2)(d\dot{q}^1-d\dot{q}^2)+
\int L_3(\dot{q}^2)\dot{q}^2d\dot{q}^2+\nonumber\\
&&\int\int
\widetilde{L}'_2(\dot{q})(\dot{q}^1-\dot{q}^2)d\dot{q}^1d\dot{q}^2+
\int\int \widetilde{L}_2(\dot{q})d\dot{q}^1d\dot{q}^2 -\nonumber\\
&& \int C_2(-q^1+f^5(t)+const)dq^1 +v_{s_{V_1}}(t)q^2\ea
where $\widetilde{L}^{'}_{2}:=\f 
{d\widetilde{L}_2(\eta)}{d\eta^l}$.$\quad \blacksquare$


\begin{thebibliography}{20}

\bibitem{marsden} R.Abraham and J.Marsden. {\it Foundations of Mechanics}. Benjamin/Cummings Publishing Company, Inc.,  1978.


\bibitem{bern} G. M. Bernstein and M.A. Lieberman. A method for obtaining a canonical Hamiltonian for nonlinear LC circuits. {\it 
IEEE Trans. Circuits and Systems}, {\bf 36} (3), pages 411-420, 1989.

 \bibitem {bir}G. D. Birkhoff. {\it Dynamical Systems}. American Mathematical Society Colloquium Publications, vol. IX, New York, 
1927.


\bibitem{bloch} A.M. Bloch and P.E. Crouch. Representations of Dirac structures on vector spaces and nonlinear LC circuits. {\it 
Differential Geometry and Control}, volume
{\bf 64} of {\it Proceedings of Symposia in Pure Mathematics}, pages 103-117. American Mathematical Society, 1999.

\bibitem{moser} R.K. Brayton and  J.K. Moser.
A theory of nonlinear networks I, II.
{\it Quarterly of Applied Mathematics}, {\bf 22}, pages 1-33,  81-104,  1964.

\bibitem{chua} L.O. Chua and J.D. McPherson. Explicit topological formulation of Lagrangian and Hamiltonian equations for nonlinear 
networks. {\it IEEE Trans. Circuits and  Systems}, {\bf 21} (2), pages 277-286, 1974.

\bibitem{graph} L.R. Foulds. {\it Graph Theory Applications}.
Springer-Verlag New York, Inc., 1992.


\bibitem{oliva} M.H. Kobayashi and W.M. Oliva. On the Birkhoff approach to classical mechanics. {\it Resenhas IME-USP} {\bf 6} (1), 
pages 1-71 , 2003.



\bibitem{maschke}  B.M. Maschke, A.J. van der Schaft and P.C. Breedveld. An intrinsic Hamiltonian formulation of the dynamics of LC 
circuits. {\it IEEE Trans. Circuits and Systems}, {\bf 42} (2), pages 73-82, 1995.
\bibitem{maschke2} B.M. Maschke and A.J. van der Schaft. The Hamiltonian formulation of energy conserving physical systems with 
external ports. {\it Archiv f\"{u}r Elektronik und Ubertragungstechnik}, {\bf 49}, pages 362-371, 1995.

\bibitem{saunders}
D.J. Saunders. {\it The Geometry of Jet Bundles}. London Mathematical Society Lecture Note Series, 142, Cambridge University Press, 
1989.

\bibitem{schaft}  A.J. van der Schaft. Implicit Hamiltonian systems with symmetry. {\it Rep. Math. Phys.}, {\bf 41}, pages 203-221, 
1998.
\end{thebibliography}
\end{document}